\documentclass[review,onefignum,onetabnum]{siamart171218}
\usepackage{amssymb}
\usepackage{amsmath}
\usepackage{xcolor}

\usepackage{amsfonts}
\usepackage{graphicx}
\usepackage{float}
\usepackage{subcaption}
\usepackage{epstopdf}
\usepackage{algorithmic}
\ifpdf
  \DeclareGraphicsExtensions{.eps,.pdf,.png,.jpg}
\else
  \DeclareGraphicsExtensions{.eps}
\fi

\newsiamremark{remark}{Remark}
\newsiamremark{hypothesis}{Hypothesis}
\crefname{hypothesis}{Hypothesis}{Hypotheses}
\newsiamthm{claim}{Claim}

\headers{Convexification for an inverse scattering problem}{V. A. Khoa, M. V. Klibanov, and L. H. Nguyen}

\title{Convexification for a 3D inverse scattering problem with the moving point source\thanks{Submitted to the editors DATE.
\funding{This work was supported by US Army Research Laboratory and US Army Research Office grant W911NF-19-1-0044. The work of V.A.K was also partially supported by the Research Foundation-Flanders (FWO) under the project named ``Approximations for forward and inverse reaction-diffusion problems related to cancer models".}}}

\author{Vo Anh Khoa\footnotemark[3]
	\and
	Michael Victor Klibanov\footnotemark[3]\;\textsuperscript{\;,}\thanks{Corresponding author.}
	\and
	Loc Hoang Nguyen\thanks{Department of Mathematics and Statistics, University of North Carolina at Charlotte, Charlotte, North Carolina 28223, USA.
	\newline \indent\: \textit{Email addresses}: \email{vakhoa.hcmus@gmail.com} (V. A. Khoa), \email{mklibanv@uncc.edu} (M. V. Klibanov), \email{loc.nguyen@uncc.edu} (L. H. Nguyen).}
}

\usepackage{amsopn}

\input{tcilatex}
\begin{document}

\maketitle

\begin{abstract}
For the first time, we develop in this paper the globally convergent
convexification numerical method for a Coefficient Inverse Problem for the
3D Helmholtz equation for the case when the backscattering data are
generated by a point source running along an interval of a straight line and
the wavenumber is fixed. Thus, by varying the wavenumber, one can
reconstruct the dielectric constant depending not only on spatial variables
but the wavenumber (i.e. frequency) as well. Our approach relies on a new
derivation of a boundary value problem for a system of coupled quasilinear
elliptic partial differential equations. This is done via an application of
a special truncated Fourier-like method. First, we prove the Lipschitz
stability estimate for this problem via a Carleman estimate. Next, using the
Carleman Weight Function generated by that estimate, we construct a globally
strictly convex cost functional and prove the global convergence to the
exact solution of the gradient projection method. Finally, our theoretical
finding is verified via several numerical tests with computationally
simulated data. These tests demonstrate that we can accurately recover all
three important components of targets of interest: locations, shapes and
dielectric constants. In particular, large target/background contrasts in
dielectric constants (up to 10:1) can be accurately calculated.
\end{abstract}




\QTP{Body Math}







\section{Introduction}

We consider a 3D Coefficient Inverse Problem (CIP) for the Helmholtz
equation in the case when the wavenumber (i.e. frequency) is fixed and the
backscattering boundary data for the inversion are generated by the point
source moving along an interval of a straight line. We develop analytically
and test computationally the so-called \emph{convexification} globally
convergent numerical method for this CIP. We call a numerical method for a
CIP \emph{globally convergent} if there is a theorem, which claims that this
method delivers at least one point in a sufficiently small neighborhood of
the correct solution without any advanced knowledge of this neighborhood. In
other words, a good first guess is not required.

The coefficient of the Helmholtz equation, i.e. the spatially distributed
dielectric constant, is the subject of the solution of our CIP. In
principle, the case when the source is moving and the frequency is fixed
enables one to consider a physically realistic problem when the dielectric
constant depends not only on spatial variables but on the frequency as well.
Indeed, if we repeat those measurements for an interval of frequencies, then
we can find the dependence of the dielectric constant on both spatial
variables and the frequency. But we assume here that the dielectric constant
depends only on spatial variables and not on the frequency. The analytical
part of this paper is devoted to the derivation of the method and its
convergence analysis. In the numerical part, we demonstrate the numerical
performance of our technique for the case of imaging of dielectric constants
of targets, which mimic antipersonnel land mines and improvised explosive
devices (IEDs).

Unlike this paper, previously the convexification method was constructed for
some CIPs for the Helmholtz equation only for the case of a single direction
of the incident plane wave with the wavenumber running over a certain
interval \cite{CAMWA,convexper,convIPnew}. We demonstrate in our numerical
studies below that in the moving source case, the convexification method
accurately images all three components of targets of interest: locations,
shapes and the target/background contrasts in the dielectric constant. This
is unlike the above mentioned previously studied case of a single direction
of the incident plane wave, which ensured only first and third components,
while shapes were not accurately imaged.

One of strengths of the convexification is that it works only with the non
overdetermined data. This means that the number $m$ of independent variables
in the data equals the number $n$ of independent variables in the unknown
coefficient, $m=n$. In particular, in our CIP $m=n=3$. On the other hand,
there are some globally convergent numerical methods for CIPs, which work
with the case $m>n.$ In this regard we refer to, e.g., \cite%
{deHoop,Kab1,Kab2}.

From the applied standpoint, we are oriented towards the problem of the
detection and identification of antipersonnel land mines and IEDs.
Reconstructions of dielectric constants from experimentally collected
backscattering data for targets mimicking these explosive devices and buried
under the ground were studied in \cite{convexper,Nguyen2018}, where a single
direction of the incident plane wave was used. Even though we do not
consider here the case of buried targets, explosives can often be located in
the air, and we model this case. In this regard we refer to the work \cite%
{KBSNF}, which analyses the experimental data collected in the field from
explosive-like targets by engineers of the US Army Research Laboratory. Some
targets in \cite{KBSNF} are located in air and some are buried in the
ground. Both here and in \cite{KBSNF} the spatially dependent dielectric
constants are subject to the solutions of CIPs. It was stated in \cite{KBSNF}
(page 33) that even though the knowledge of the dielectric constant alone is
insufficient to identify an explosive, one can still hope that this
knowledge might serve as an important piece of information, additional to
the conventional ones, to help better identify explosives and thus, to
decrease the false alarm rate.

Any CIP is both nonlinear and ill-posed. Therefore, a conventional least
squares cost functional for this problem is, as a rule, non convex; see,
e.g. \cite{Chavent,Gonch} for some works in which least squares cost
functionals are applied to solve CIPs. The non convexity, combined with the
ill-posedness, causes the presence of many local minima and ravines; see,
e.g., \cite{Scales} for a convincing numerical example of multiple local
minima. Since a minimization procedure can stop at any local minimum, there
is no guarantee that the solution obtained via the minimization process
applied to that functional is indeed close to the correct one. The only
exception is sometimes the case when the starting point of iterations is
located in a sufficiently small neighborhood of the correct solution. We
call the latter \emph{local convergence}. However, a good first guess about
the solution is rarely available in applications.

The above motivates this research group to work on the convexification
approach. The roots of the convexification are in the method of Carleman
estimates for CIPs. This method was originated in the work \cite{BukhKlib}.
The idea of \cite{BukhKlib} led to many publications of many authors. Since
this paper is not a survey of the method of \cite{BukhKlib}, we refer for
brevity only to the books \cite{BK,BY,KT} and the survey \cite{Ksurvey}. We
also note that initially the method of \cite{BukhKlib} was created
exclusively for proofs of uniqueness and stability theorems for CIPs.

In the convexification, one constructs a weighted Tikhonov-like functional.
The weight is the Carleman Weight Function, i.e. the function involved as
the weight in the Carleman estimate for the corresponding PDE operator.
Given a convex bounded set $D\left( \beta \right) \subset H$ of an arbitrary
diameter $\beta >0$ in an appropriate Hilbert space $H$, one can choose the
parameter $\lambda >0$ of the Carleman Weight Function such that the strict
convexity of the functional on that set is ensured. Thus, the phenomenon of
local minima does not occur. Furthermore, starting from the publication \cite%
{Bak}, all works on the convexification include theorems, which claim
convergence of the gradient projection method of the minimization of that
functional to the correct solution of the corresponding CIP if starting from
an arbitrary point of $D\left( \beta \right) $. Given that the diameter $%
\beta >0$ of $D\left( \beta \right) $ is an arbitrary one, this is \emph{%
global convergence}.

Initial publications on the convexification \cite{BKconv,Klib97} were only
theoretical ones. More recently, however, the work \cite{Bak} has clarified
some points, which were preventing one from numerical studies. As a result,
the most recent works \cite{CAMWA,convexper,convIPnew,EIT,timedomain}, so as
the current one, contain both analytical and numerical studies of the
convexification. In particular, an accurate performance of the
convexification on experimental backscattering data was demonstrated in \cite%
{convexper,convIPnew}. One of important conclusions of these numerical
studies is that even though the theory requires large values of the
parameter $\lambda ,$ accurate numerical results can be obtained for
reasonable values $\lambda \in \left[ 1,3\right] .$

We also refer to publications \cite{Baud,Baud1}, where a different
version of the convexification is used to develop globally convergent
numerical methods for CIPs for some hyperbolic PDEs. Just as in the
convexification, the Carleman Weight Functions play a pivotal role in \cite{Baud,Baud1}.

As to the CIPs with the fixed wavenumber, we refer to
numerical procedures developed during a long standing effort by the group of
Novikov since about 1988 \cite{N1}; also, see, e.g. \cite{Ag,Al,N2,N3}. The
statements of CIPs in these publications are different from ours. These
reconstruction techniques are also different from the convexification. By
the above definition, they can also be called ``globally convergent
numerical methods", since neither of them requires a good first guess for
the solution while still delivering the exact coefficient in the end. An
interesting feature of \cite{Ag,N3} is that these publications consider the
case of the non overdetermined data for the reconstruction of the potential
of the Schr\"odinger equation at high values of the wavenumber. Another
noticeable feature of \cite{Ag} is that the data there are phaseless.
Corresponding numerical results can be found in \cite{Ag,Al}.

The main new elements of this paper are:

\begin{enumerate}
\item It is the first time when the convexification is applied to a CIP for
the Helmholtz equation in the case when the data are generated by the moving
point source and the wavenumber is fixed.

\item We prove the Lipschitz stability estimate for an overdetermined
boundary value problem for an auxiliary system of coupled quasilinear
elliptic PDEs. This \ result is interesting in its own right.

\item In the proof of the central theorem about the global strict convexity
of our above mentioned weighted Tikhonov-like functional we do not subtract
the boundary data from the solution of that system of quasilinear elliptic
PDEs. In other words, we do not arrange boundary conditions for that
difference to be equal to zero.

\item We prove the Lipschitz stability of minimizers of our weighted
Tikhonov-like functional with respect to small perturbations of the data as
well as \textquotedblleft Lipschitz-like" convergence rate of the gradient
projection method (the latter converges globally). These results are
stronger than those of all previous works on the convexification \cite%
{CAMWA,convexper,convIPnew,EIT,timedomain}, where the weaker H\"{o}lder
stability of minimizers and the \textquotedblleft H\"{o}lder-like"
convergence rates were proven.

\item The numerical results are new.
\end{enumerate}

We now discuss one important issue related to this paper. We refer
to Remarks 3.1 for a further discussion of this issue. It is well known that
the CIP under consideration is a \textbf{\emph{substantially challenging one}}. 
This is the fundamental underlying reason why we actually replace
our CIP with an approximate one, which we call \textquotedblleft approximate
mathematical model". Indeed, even the uniqueness of our CIP can be proven,
at least at this particular moment of time, only within the framework of
this model: it follows immediately from Theorem 3.2, also; see \cref{sec:2}. In
fact, the major part of our analytical treatment works \textbf{\emph{only within the
framework of our approximate mathematical model}}. This model amounts to the
truncation of a certain Fourier series. However, due to both the
ill-posedness and the nonlinearity of our CIP, we cannot provide an
analytical estimate as $N\rightarrow \infty $ of the accuracy of
the solution of the CIP which results from that approximation. Here, $N$ 
 is the number of terms of that truncated Fourier series.
Nevertheless, we analyze numerically the $N-$dependence of the
approximation accuracy of a function generated by the solution of the
forward problem. Based on this analysis, we choose the number $N$
numerically and consider this choice as an optimal one, see \cref{sec:experiments}
and \cref{fig:0a}. 

In our numerical studies (\cref{sec:alg,sub:7.2.3}), we generate the
data for our CIP via the solution of the forward problem without any
truncation of any Fourier series. Next, we apply our method to these 
data to find those $N$ Fourier coefficients
\textquotedblleft pretending" that these data are governed by our
approximate mathematical model. We believe that a good accuracy of our
reconstruction results justifies our model from the numerical standpoint. 

We point out that a similar situation is quite \textbf{typical} in the field
of Inverse Problems: when a very challenging problem is replaced with an
approximate mathematical model for the sake of a numerical method.
Corresponding numerical results are usually good ones. In this regard, we
refer to some works of other authors \cite{GN,Kab1,Kab2} as well as to our
previous works on the convexification \cite{CAMWA,convexper,convIPnew,EIT,timedomain}. In addition, the above mentioned works \cite{Ag,Al,N2,N3} actually also use certain truncations of Fourier type transforms at some steps, i.e., just as ourselves, they use some approximate mathematical models. The computational experience of all these cited works, including the
current paper, tells one that the number $N$, including its
analogs, can be successfully chosen numerically.

At the start of our work on this CIP, we had two choices: either not
to work on it, due to the above mentioned substantial challenges, or to
figure out an approximate mathematical model and then to develop a globally
convergent numerical method for the latter. Keeping in mind the importance
of the above discussed application, we have chosen the second option. 

The paper is organized as follows. In the next section, we detail the CIP we
work with. Then, in \cref{sec:3} an auxiliary boundary value problem is
derived for a system of coupled quasilinear elliptic PDEs and the Lipschitz
stability of this problem is proven. Next, in \cref{sec:4} the above
mentioned weighted Tikhonov-like functional is constructed. The central
theorem about the global strict convexity of this functional is proven in %
\cref{sec:5}. In addition, the Lipschitz stability of minimizers is also
proven in \cref{sec:5}. In \cref{sec:6} we establish the global convergence
of the gradient projection method and present its convergence rate.
Numerical studies are described in \cref{sec:7}. 

\section{Statement of the Coefficient Inverse Problem}

\label{sec:2} We model the propagation of the electric wave field by the
Helmholtz equation instead of the Maxwell's equations. This modeling was
numerically justified in Appendix of the paper \cite{KNN}. Such a
mathematical model is true at least for rather simple medium consisting of a
homogeneous background and a few embedded inclusions. Besides, good
accuracies of reconstructions obtained by our research group from
experimental data in publications \cite{convexper,JCP,Nguyen2018}, where the
Helmholtz equation was used to model the wave propagation process, speak in
favor of this modeling.

Denote $\mathbf{x}=(x,y,z)\in \mathbb{R}^{3}$. Let the number $R>0.$ We
define the cube $\Omega \subset \mathbb{R}^{3}$ as 
\begin{equation}
\Omega =\left\{ \mathbf{x}:\left\vert x\right\vert ,\left\vert y\right\vert
,\left\vert z\right\vert <R\right\} .  \label{200}
\end{equation}%
Let $\Gamma \subset \partial \Omega $ be the lower part of the boundary of $%
\Omega $ where measurements of the backscatter data are conducted, 
\begin{equation}
\Gamma :=\left\{ \mathbf{x}:\left\vert x\right\vert ,\left\vert y\right\vert
<R,z=-R\right\} .  \label{201}
\end{equation}%
Let $c:=c(\mathbf{x})\in \lbrack 1,\infty )$ be a sufficiently smooth
function that represents the dielectric function of the medium. We assume
that 
\begin{equation}
\begin{cases}
c\left( \mathbf{x}\right) \geq 1 & \text{in }\mathbb{R}^{3}, \\ 
c\left( \mathbf{x}\right) =1 & \text{in }\mathbb{R}^{3}\backslash \Omega .%
\end{cases}
\label{2020}
\end{equation}%
Here, $k>0$ is the wavenumber. The function $c$ is the spatially distributed
and $k$dependent dielectric constant. The second assumption in \eqref{2020}
means that we have vacuum outside of the domain of interest $\Omega $.

Let $a$ and $d$ be two numbers such that $d>R$ and $a>0$. We define the 
\emph{line of sources }as 
\begin{equation}
L_{\text{src}}:=\left\{ \left( \alpha ,0,-d\right) :-a\leq \alpha \leq
a\right\} .  \label{2001}
\end{equation}%
Obviously, this line is parallel to the $x$--axis. The distance from $L_{%
\mathbf{src}}$ to $\Gamma $ is $d$, and the length of our the line of
sources is $2a$. Since $d>R$, then $L_{\text{src}}\cap \overline{\Omega }%
=\varnothing $. Thus, for each $\alpha \in \lbrack -a,a]$ the corresponding
point source is $\mathbf{x}_{\alpha }:=(\alpha ,0,-d)\in L_{\text{src}}$.

First, we formulate the forward problem. Let $k=const.>0$ and assume that
the function $c$ is known. For each source position $\mathbf{x}_{\alpha }\in
L_{\text{src}}$ the forward problem is: 
\begin{align}
& \Delta u+k^{2}c\left( \mathbf{x}\right) u=-\delta \left( \mathbf{x}-%
\mathbf{x}_{\alpha }\right) \quad \text{in }\mathbb{R}^{3},
\label{eq:forward1} \\
& \lim_{r\rightarrow \infty }r\left( \partial _{r}u-\text{i}ku\right)
=0\quad \text{for }r=\left\vert \mathbf{x}-\mathbf{x}_{\alpha }\right\vert ,%
\text{i}=\sqrt{-1}.  \label{eq:forward2}
\end{align}

Conditions \eqref{eq:forward1}--\eqref{eq:forward2} form the Helmholtz
equation with the Sommerfeld radiation condition at the infinity. Let $u_{0}(%
\mathbf{x},\alpha )$ be the solution of \eqref{eq:forward1}--%
\eqref{eq:forward2} with $c\equiv 1$, 
\begin{equation}
u_{0}\left( \mathbf{x},\alpha \right) =\frac{\exp \left( \text{i}k\left\vert 
\mathbf{x}-\mathbf{x}_{\alpha }\right\vert \right) }{4\pi \left\vert \mathbf{%
x}-\mathbf{x}_{\alpha }\right\vert }.  \label{1}
\end{equation}%
Using the Helmholtz equation for $u_{0,\alpha }=u_{0}(\mathbf{x},\alpha )$,
we obtain from \eqref{eq:forward1}--\eqref{eq:forward2} 
\begin{align*}
& \Delta \left( u-u_{0,\alpha }\right) +k^{2}\left( u-u_{0,\alpha }\right)
=-k^{2}\left( c\left( \mathbf{x},k\right) -1\right) u\quad \text{in }\mathbb{%
R}^{3}, \\
& \lim_{r\rightarrow \infty }r\left[ \partial _{r}\left( u-u_{0,\alpha
}\right) -\text{i}k\left( u-u_{0,\alpha }\right) \right] =0\quad \text{for }%
r=\left\vert \mathbf{x}-\mathbf{x}_{\alpha }\right\vert .
\end{align*}%
In view of the fact that $c(\mathbf{x})=1$ in $\mathbb{R}^{3}\diagdown
\Omega $, we thus find that the solution $u$ to the system %
\eqref{eq:forward1}--\eqref{eq:forward2} satisfies the so-called
Lippmann--Schwinger equation (see, e.g. \cite[Section 8.2]{Colton1992}),
which reads for all $\mathbf{x}\in \mathbb{R}^{3}$ as 
\begin{equation}
u\left( \mathbf{x},\alpha \right) =u_{0}\left( \mathbf{x},\alpha \right)
+k^{2}\int_{\Omega }\frac{\exp \left( \text{i}k\left\vert \mathbf{x}-\mathbf{%
x}^{\prime }\right\vert \right) }{4\pi \left\vert \mathbf{x}-\mathbf{x}%
^{\prime }\right\vert }\left( c\left( \mathbf{x}^{\prime }\right) -1\right)
u\left( \mathbf{x}^{\prime },\alpha \right) d\mathbf{x}^{\prime }.
\label{302}
\end{equation}

\begin{figure}[htbp]
\begin{center}
\includegraphics[width = 0.48\textwidth]{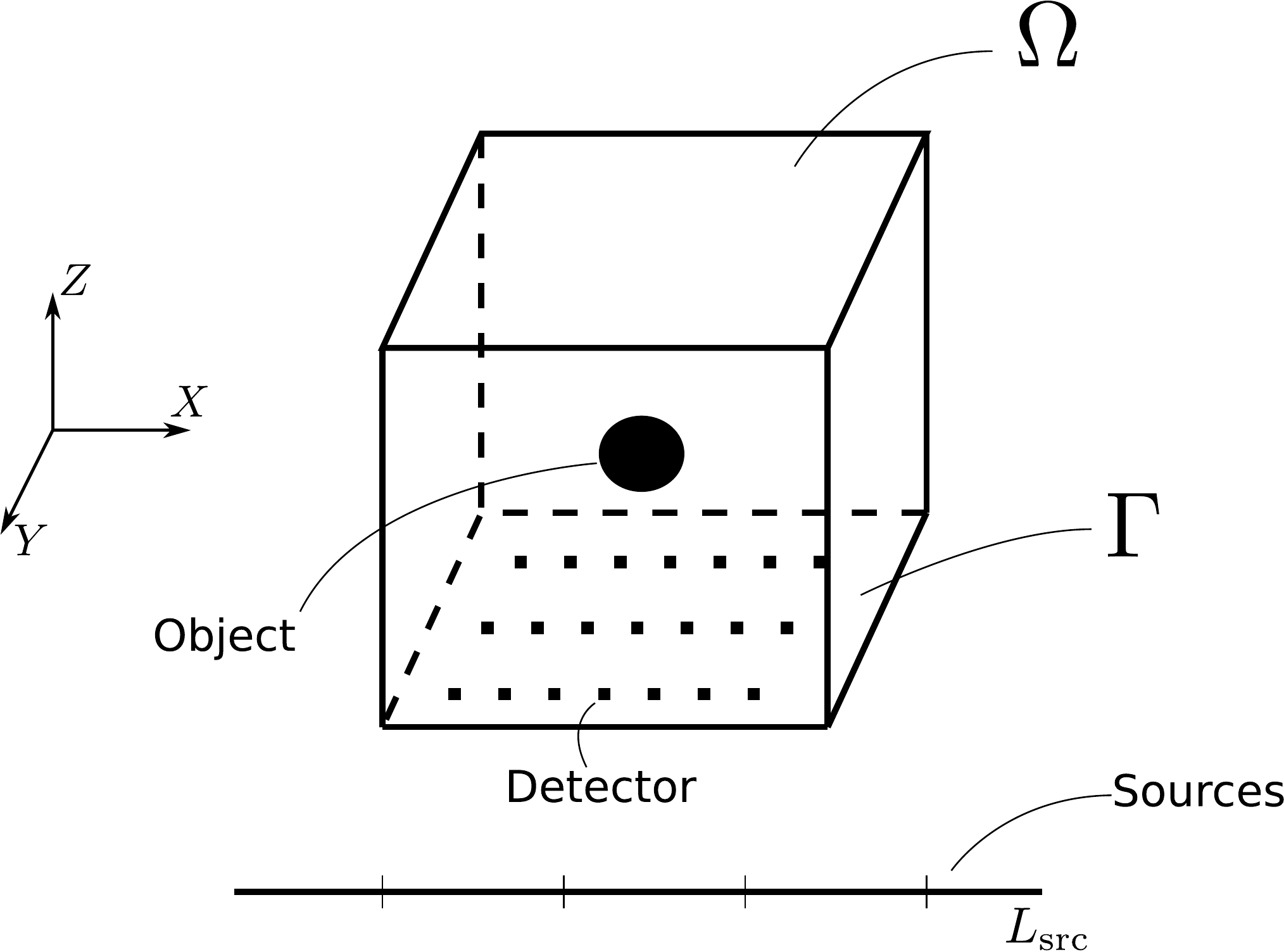}
\end{center}
\caption{A schematic diagram of data collection for our CIP. The wave field
is generated by point sources $\mathbf{x}_{\protect\alpha}\in L_{\text{src}}$%
. We measure the backscattering wave field at an array of detectors on the
lower side $\Gamma$ of the cube $\Omega$.}
\label{fig:1}
\end{figure}

We now pose the CIP which we solve in this paper. The schematic diagram of
measurements for this problem is illustrated in \cref{fig:1}. 


\textbf{Coefficient Inverse Problem (CIP). }\emph{Given }$k>0$\emph{,
determine the coefficient }$c(\mathbf{x})$\emph{\ for }$\mathbf{x}\in \Omega 
$\emph{\ in the system \eqref{eq:forward1}--\eqref{eq:forward2}, assuming
that the following function }$F(\mathbf{x},\mathbf{x}_{\alpha })$ \emph{is
given }%
\begin{equation}
F(\mathbf{x},\mathbf{x}_{\alpha })=u(\mathbf{x},\alpha )\quad \text{for }%
\mathbf{x}\in \Gamma ,\mathbf{x}_{\alpha }\in L_{\text{src}},  \label{303}
\end{equation}%
\emph{where }$u(\mathbf{x},\alpha )$\emph{\ is the solution to %
\eqref{eq:forward1}--\eqref{eq:forward2}}.

Physically, to reconstruct the dielectric function $c$ of objects in $\Omega 
$, one sends the incident wave field from the source $\mathbf{x}_{\alpha }$.
This wave scatters when hitting the objects. Then, one measures the
backscattering wave on the square $\Gamma $. And the data \eqref{303} are
used to reconstruct the unknown dielectric constant inside the cube $\Omega $%
.

Uniqueness of this CIP\ is a long standing open problem. Currently
uniqueness can be proven by the method of \cite{BukhKlib} only if the right
hand side of equation \eqref{eq:forward1} is not vanishing in $\overline{%
\Omega }.$ Nevertheless, uniqueness within the framework of our approximate
mathematical model (Remarks 3.1) follows immediately from Theorem 3.2.

\textbf{Remarks 2.1}.

1. \emph{In this work, we are not interested in a specification of
smoothness condition imposed on the function $c(\mathbf{x}).$ Thus, $c(%
\mathbf{x})$ is supposed to be sufficiently smooth with respect to $\mathbf{x%
}$. Some particular discussions concerning this matter can be found in,
e.g., \cite{CAMWA} and references therein, where the smoothness of $c$ is
essential for the asymptotic behavior of the solution $u$ to the forward
problem \eqref{eq:forward1}--\eqref{eq:forward2}. We also note that in
studies of CIPs the smoothness conditions are usually not of a considerable
concern, see e.g. \cite[Theorem 4.1]{Rom}.}

2. \emph{To solve the forward problem \eqref{eq:forward1}--%
\eqref{eq:forward2} using the integral equation \eqref{302} for all $\mathbf{%
x}\in \Omega $, we rely on numerical methods commenced in \cite{Vainikko2000}%
. This way enables us to extract information of $\left. u\left( \mathbf{x}%
,\alpha \right) \right\vert _{\Gamma }$, and by repeating this process for
each $\alpha \in \lbrack -a,a]$ we obtain computationally the simulated data %
\eqref{303}.}

\section{An Auxiliary System of Coupled Quasilinear Elliptic Equations}

\label{sec:3}

\subsection{An equation without the unknown coefficient}

\label{sec:3.1}

Observe that since $L_{\text{src}}$ is located outside of $\overline{\Omega }
$, then the point source $\mathbf{x}_{\alpha }=(\alpha ,0,-d)$ is not in $%
\overline{\Omega }$. Hence, \eqref{eq:forward1}--\eqref{eq:forward2} imply
that for each $\alpha \in \lbrack -a,a]$ 
\begin{equation}
\Delta u+k^{2}c\left( \mathbf{x}\right) u=0\quad \text{in }\Omega .
\label{eq:aux}
\end{equation}%
We now define the function $\log u\left( \mathbf{x},\alpha \right)$
. The conformal Riemannian metric generated by the function $
c\left( \mathbf{x}\right)$ is %
\begin{equation*}
d\tau =\sqrt{c\left( \mathbf{x}\right) }\left\vert d\mathbf{x}\right\vert ,%
\text{ }\left\vert d\mathbf{x}\right\vert =\sqrt{\left( dx\right)
^{2}+\left( dy\right) ^{2}+\left( dz\right) ^{2}}.
\end{equation*}%
Following \cite{KR}, we assume that geodesic lines generated by this
metric and originated at sources $\mathbf{x}_{\alpha }\in L_{\text{src}}$  are regular. In other words, for each point $\mathbf{x}\in 
\mathbb{R}^{3}$ and for each point $\mathbf{x}_{\alpha }\in L_{\text{src}}$ there exists unique geodesic line $\Gamma \left( 
\mathbf{x},\mathbf{x}_{\alpha }\right) $ connecting them. A
sufficient condition of the regularity of geodesic lines was established in \cite{Rom1} and can also be found in \cite{KR}. The
travel time along $\Gamma \left( \mathbf{x},\mathbf{x}_{\alpha }\right) $  is
\begin{equation*}
\tau \left( \mathbf{x},\mathbf{x}_{\alpha }\right) =\dint\limits_{\Gamma
\left( \mathbf{x},\mathbf{x}_{\alpha }\right) }\sqrt{c\left( \xi %
\right) }ds.
\end{equation*}%
Tentatively, we denote $u=u\left( \mathbf{x},k,\alpha \right) .$ 
It was established in \cite{KR} that, under certain conditions imposed on $c\left( \mathbf{x}\right) ,$ which we do not discuss here (Remark
2.1), asymptotic behavior of the function $u\left( \mathbf{x},k,\alpha
\right) $ at $k\rightarrow \infty $ is:
\begin{align}
 u\left( \mathbf{x},k,\alpha \right) =A\left( \mathbf{x},\alpha \right) \exp
\left( \text{i}k\tau \left( \mathbf{x},\mathbf{x}_{\alpha }\right) \right) \left[
1+\mathcal{O}\left( 1/k\right) \right] \quad \text{ for all } \left( \mathbf{x},\alpha \right) \in \overline{\Omega }\times \left[ -a,a\right] ,  \label{2000} 
\end{align}where the function $A\left( \mathbf{x},\alpha \right) >0.$ 
 Let $\overline{k}>1$ be a number. Assuming that $\overline{k}$  is sufficiently large, and that $k\geq \overline{k},$ we
obtain from (\ref{2000}) that $u\left( \mathbf{x},k,\alpha
\right) \neq 0$ for all $\left( \mathbf{x},k,\alpha \right) \in 
\overline{\Omega }\times \left[ \overline{k},\infty \right) \times \left[
-a,a\right] .$ Denoting the term $\mathcal{O}\left( 1/k\right) $ in
(\ref{2000}) as $\mathcal{O}\left( 1/k\right) =s\left( \mathbf{x},k,\alpha \right) $, we naturally assume that $\left\vert s\left( \mathbf{x},k,\alpha
\right) \right\vert <1.$ Hence, using (\ref{2000}), we uniquely
define the function $\log u\left( \mathbf{x},k,\alpha \right) $ as\begin{equation}
\log u\left( \mathbf{x},k,\alpha \right) =\text{i}k\tau \left( \mathbf{x},\mathbf{x}_{\alpha }\right) +\ln A\left( \mathbf{x},\alpha \right)
+\dsum\limits_{n=1}^{\infty }\frac{\left( -1\right) ^{n-1}}{n}s^{n}\left( 
\mathbf{x},k,\alpha \right) ,  \label{2002}
\end{equation}for all $\left( \mathbf{x},k,\alpha \right) \in \overline{\Omega }\times \left[ \overline{k},\infty \right) \times \left[ -a,a\right] .$  In (\ref{2000}) the series is obviously taken from the power
series expansion of the function $\log \left( 1+s\left( \mathbf{x},k,\alpha
\right) \right) .$

Now, suppose that $k$ is not so large, but still $k\geq 
\underline{k}>0$  where the number $\underline{k}$ is such
that %
\begin{equation}
u\left( \mathbf{x},k,\alpha \right) \neq 0 \quad \text{for all } \left( 
\mathbf{x},k,\alpha \right) \in \overline{\Omega }\times \left[ \underline{k}%
,\infty \right) \times \left[ -a,a\right] .  \label{2003}
\end{equation}%
 Then, using an idea of \cite{KlibNg}, we define the
function $\varphi \left( \mathbf{x},k,\alpha \right) $ as%
\begin{equation}
\varphi \left( \mathbf{x},k,\alpha \right) =-\dint\limits_{k}^{\overline{k}}%
\frac{\partial _{k}u\left( \mathbf{x},\eta ,\alpha \right) }{u\left( \mathbf{%
x},\eta ,\alpha \right) }d\eta +\log u\left( \mathbf{x},\overline{k},\alpha
\right) ,\text{ } \forall \left( \mathbf{x},k,\alpha \right) \in \overline{%
\Omega }\times \left[ \underline{k},\infty \right) \times \left[ -a,a\right]
,  \label{2004}
\end{equation}%
where $\log u\left( \mathbf{x},\overline{k},\alpha \right) $\ is defined in (\ref{2002}). Hence, $\left( \partial _{k}u-u\partial
_{k}\varphi \right) \left( \mathbf{x},k,\alpha \right) =0.$ 
Multiplying both sides of the latter by $e^{-\varphi },$ we obtain 
$\partial _{k}\left( ue^{-\varphi }\right) \left( \mathbf{x},k,\alpha
\right) =0.$ Since $\varphi \left( \mathbf{x},\overline{k},\alpha
\right) =u\left( \mathbf{x},\overline{k},\alpha \right) ,$ then $u\left( \mathbf{x},k,\alpha \right) =\exp \left( \varphi \left( \mathbf{x},k,\alpha \right) \right) .$ This uniquely defines the function $\log u\left( \mathbf{x},k,\alpha \right) =\varphi \left( \mathbf{x},k,\alpha
\right) $ as long as (\ref{2003}) holds. Finally, we note that we
use below only derivatives of the function $u\left( \mathbf{x},k,\alpha
\right) ,$ which means that we do not use $\log u$
\textquotedblleft directly".

In all our above cited previous publications about numerical methods for
CIPs for the Helmholtz equation we have not observed numerically such values
of the function $\left\vert u\left( \mathbf{x},k,\alpha \right) \right\vert $
which would be close to zero. The same is true for the current paper. Thus,
we assume below that the fixed number $k$ we work with is such $k\in \left[ 
\underline{k},\infty \right) .$ Hence, by (\ref{2003})--(\ref{2004}),
the function $\log u\left( \mathbf{x},k,\alpha \right) =\varphi \left( 
\mathbf{x},k,\alpha \right) $ is uniquely defined. Thus, we assume
below that
\begin{equation*}
u\left( \mathbf{x},\alpha \right) \neq 0 \quad \text{for all } \left( \mathbf{x},\alpha
\right) \in \overline{\Omega }\times \left[ -a,a\right] .
\end{equation*}%
We set 
\begin{equation}
\log u_{0}\left( \mathbf{x},\alpha \right) =\text{i}k\left\vert \mathbf{x}-%
\mathbf{x}_{\alpha }\right\vert -\ln \left( 4\pi \left\vert \mathbf{x}-%
\mathbf{x}_{\alpha }\right\vert \right) .  \label{3.2}
\end{equation}

Denote $v_{0}(\mathbf{x},\alpha )=u(\mathbf{x},\alpha )/u_{0}(\mathbf{x}%
,\alpha )$ and define the function $v(\mathbf{x},\alpha )$ as 
\begin{equation}
v(\mathbf{x},\alpha )=\log v_{0}(\mathbf{x},\alpha )=\log u(\mathbf{x}%
,\alpha )-\log u_{0}(\mathbf{x},\alpha )\quad \text{for }\mathbf{x}\in
\Omega ,\alpha \in \lbrack -a,a].  \label{300}
\end{equation}%
Obviously, 
\begin{equation}
\nabla v\left( \mathbf{x},\alpha \right) =\frac{\nabla v_{0}\left( \mathbf{x}%
,\alpha \right) }{v_{0}\left( \mathbf{x},\alpha \right) },\text{ }\Delta
v\left( \mathbf{x},\alpha \right) =\frac{\Delta v_{0}\left( \mathbf{x}%
,\alpha \right) }{v_{0}\left( \mathbf{x},\alpha \right) }-\left( \frac{%
\nabla v_{0}\left( \mathbf{x},\alpha \right) }{v_{0}\left( \mathbf{x},\alpha
\right) }\right) ^{2}.  \label{3.02}
\end{equation}%
Using (\ref{3.02}), we obtain the equation for $v$: 
\begin{equation}
\Delta v+\left( \nabla v\right) ^{2}+2\nabla v\cdot \nabla \left( \log u_{0}(%
\mathbf{x},\alpha )\right) =-k^{2}(c\left( \mathbf{x},k\right) -1),\text{ }%
\mathbf{x}\in \Omega .  \label{3.03}
\end{equation}%
Differentiating \eqref{3.03} with respect to $\alpha $, we obtain 
\begin{equation}
\Delta \partial _{\alpha }v+2\nabla v\cdot \nabla \partial _{\alpha
}v+2\nabla \partial _{\alpha }v\cdot \tilde{\mathbf{x}}_{\alpha }+2\widehat{%
\mathbf{x}}_{\alpha }\cdot \nabla v=0\;\text{ for all }\mathbf{x}\in \Omega .
\label{eq:aux1}
\end{equation}%
Recall that $\mathbf{x}-\mathbf{x}_{\alpha }=(x-\alpha ,y,z+d)$. We have the
following notations in (\ref{eq:aux1}): 
\begin{align}
\tilde{\mathbf{x}}_{\alpha }& =\frac{\text{i}k\left( \mathbf{x}-\mathbf{x}%
_{\alpha }\right) }{\left\vert \mathbf{x}-\mathbf{x}_{\alpha }\right\vert }-%
\frac{\mathbf{x}-\mathbf{x}_{\alpha }}{\left\vert \mathbf{x}-\mathbf{x}%
_{\alpha }\right\vert ^{2}},  \label{400} \\
\widehat{\mathbf{x}}_{\alpha }& =\frac{\text{i}k}{\left\vert \mathbf{x}-%
\mathbf{x}_{\alpha }\right\vert ^{3}}\left( -y^{2}-\left( z+d\right)
^{2},\left( x-\alpha \right) y,\left( x-\alpha \right) z\right)  \notag \\
& -\frac{1}{\left\vert \mathbf{x}-\mathbf{x}_{\alpha }\right\vert ^{4}}%
\left( \left( x-\alpha \right) ^{2}-y^{2}-\left( z+d\right) ^{2},2\left(
x-\alpha \right) y,2\left( x-\alpha \right) z\right) .  \notag
\end{align}

The notion behind this differentiation is to get rid of the $\alpha $%
-independent dielectric function $c$ in \eqref{3.03} and thus, the auxiliary
equation depends only on $v$ and $\partial _{\alpha }v$ is presented in %
\eqref{eq:aux1}. This approach is actually very similar with the first step
of the method of \cite{BukhKlib,Ksurvey}, which, however, was initially
proposed only for proofs of uniqueness theorems. To deal with the variable $%
\alpha $ in \eqref{eq:aux1}, we rely below on a special orthonormal basis
with respect to $\alpha $ to reduce \eqref{eq:aux1} to a system of coupled
elliptic quasilinear PDEs.

\subsection{A special Fourier basis}

\label{sec:3.2} To approximately solve the auxiliary problem \eqref{eq:aux1}%
, we use a truncated Fourier series. To do this, we use a special
orthonormal basis in $L^{2}(-a,a)$ denoted by $\left\{ \Psi _{n}(\alpha
)\right\} _{n=0}^{\infty }$, $\alpha \in (-a,a).$ This basis was first
constructed in \cite{Klibanov2017}.

For each $n\in \mathbb{N}$, let $\varphi _{n}(\alpha )={\alpha }%
^{n}e^{\alpha }$ for $\alpha \in \lbrack -a,a]$. Observe that the set $%
\left\{ \varphi _{n}(\alpha )\right\} _{n\in \mathbb{N}}$ is linearly
independent and complete in $L^{2}(-a,a)$. Using the classical Gram--Schmidt
orthonormalization procedure, we obtain the orthonormal basis $\left\{ \Psi
_{n}(\alpha )\right\} _{n\in \mathbb{N}}$ in $L^{2}(-a,a)$. see, e.g. This
basis possesses the following main properties \cite{Klibanov2017}:

\begin{itemize}
\item $\Psi _{n}\in C^{\infty }[-a,a]$ for all $n\in \mathbb{N}$;

\item Let $s_{mn}=\left\langle \Psi _{n}^{\prime },\Psi _{m}\right\rangle $
where $\left\langle \cdot ,\cdot \right\rangle $ denotes the scalar product
in $L^{2}(-a,a)$. Then the square matrix $S_{N}=\left( s_{mn}\right)
_{m,n=0}^{N-1}$ is invertible for any $N$ since 
\begin{equation*}
s_{mn}=%
\begin{cases}
1 & \text{if }n=m, \\ 
0 & \text{if }n<m.%
\end{cases}%
\end{equation*}
\end{itemize}

We note that neither classical orthogonal polynomials nor the classical basis
of trigonometric functions do not hold the second property. This is because
in any of these two the first column of the integer $N\geq 1$ the matrix $%
S_{N}$ would be identically zero. By virtue of this property, the matrix $%
S_{N}$ is actually an upper diagonal matrix with $\text{det}(S_{N})=1$. Hence, the
inverse matrix $S_{N}^{-1}$ exists.

Consider the auxiliary function $v(\mathbf{x},\alpha )$ that we have defined
in \cref{sec:3.1}. Given $N\geq 1$, our truncated Fourier series for $v$ is
written as 
\begin{equation}
v\left( \mathbf{x},\alpha \right) =\sum_{n=0}^{N-1}\left\langle v\left( 
\mathbf{x},\cdot \right) ,\Psi _{n}\left( \cdot \right) \right\rangle \Psi
_{n}\left( \alpha \right) \text{ for }\mathbf{x}\in \Omega ,\alpha \in \left[
-a,a\right] .  \label{eq:truncated}
\end{equation}%
Actually the sign \textquotedblleft $\approx "$ should be
used in (\ref{eq:truncated}). However, we use \textquotedblleft =" for the
further convenience of our work with our approximate mathematical model; see
Remarks 3.1 about this model.

\textbf{Remarks 3.1}.

\begin{enumerate}
\item \emph{The} \emph{representation (\ref{eq:truncated}) is an
approximation of the function }$v\left( \mathbf{x},\alpha \right) $\emph{\
since the rest of the Fourier series is not counted here. Furthermore, we
assume that the }$\alpha -$\emph{derivative }$\partial _{\alpha }v\left( 
\mathbf{x},\alpha \right) $\emph{\ can be obtained via the term-by-term
differentiation of the right-hand side of (\ref{eq:truncated}) with respect
to }$\alpha .$\emph{\ Next, we suppose that the substitution of (\ref%
{eq:truncated}) and its }$\alpha -$\emph{derivative in the left hand side of
equation (\ref{eq:aux1}) give us zero in its right hand side. In addition,
we assume that the substitution of (\ref{eq:truncated}) in the left hand
side of (\ref{3.03}) provides us with the exact coefficient }$c\left( 
\mathbf{x}\right) $\emph{\ in its right hand side. Finally, we impose in
Section 3.2 the boundary condition (\ref{301}) on }$\partial \Omega
\diagdown \Gamma $.

\item \emph{The assumptions of item 1 form our \textbf{approximate
mathematical model}. We cannot prove convergence as }$N\rightarrow \infty .$%
\emph{\ Indeed, such a result is very hard to prove due to both the
nonlinearity and the ill-posedness of our CIP. Therefore, our goal below is
to find }\emph{spatially dependent Fourier coefficients }$v_{n}\left( 
\mathbf{x}\right) =\left\langle v\left( \mathbf{x},\cdot \right) ,\Psi
_{n}\left( \cdot \right) \right\rangle $\emph{. The number }$N$\emph{\
should be chosen numerically; see \cref{sec:experiments} and \cref{fig:0a}. }

\item \emph{Everywhere below we work only within the framework of this
approximate mathematical model. As it was pointed out in Introduction, the
fundamental underlying reason why we are accepting this model is that the
original CIP is an extremely challenging one.}
\end{enumerate}

We now substitute \eqref{eq:truncated} into \eqref{eq:aux1} to get 
\begin{align*}
& \Delta \left( \sum_{n=0}^{N-1}v_{n}\left( \mathbf{x}\right) \Psi
_{n}^{\prime }\left( \alpha \right) \right) +2\nabla \left(
\sum_{n=0}^{N-1}v_{n}\left( \mathbf{x}\right) \Psi _{n}\left( \alpha \right)
\right) \cdot \nabla \left( \sum_{n=0}^{N-1}v_{n}\left( \mathbf{x}\right)
\Psi _{n}^{\prime }\left( \alpha \right) \right) \\
& +2\nabla \left( \sum_{n=0}^{N-1}v_{n}\left( \mathbf{x}\right) \Psi
_{n}^{\prime }\left( \alpha \right) \right) \cdot \tilde{\mathbf{x}}_{\alpha
}+2\bar{\mathbf{x}}_{\alpha }\cdot \nabla \left( \sum_{n=0}^{N-1}v_{n}\left( 
\mathbf{x}\right) \Psi _{n}\left( \alpha \right) \right) =0.
\end{align*}%
This equation is equivalent with: 
\begin{align}
& \sum_{n=0}^{N-1}\Psi _{n}^{\prime }\left( \alpha \right) \Delta
v_{n}\left( \mathbf{x}\right) +2\sum_{n=0}^{N-1}\sum_{l=0}^{N-1}\Psi
_{n}\left( \alpha \right) \Psi _{l}^{\prime }\left( \alpha \right) \nabla
v_{n}\left( \mathbf{x}\right) \cdot \nabla v_{l}\left( \mathbf{x}\right)
\label{eq:3.4} \\
& +2\Psi _{n}^{\prime }\left( \alpha \right) \sum_{n=0}^{N-1}\nabla
v_{n}\left( \mathbf{x}\right) \cdot \tilde{\mathbf{x}}_{\alpha }+2\Psi
_{n}\left( \alpha \right) \sum_{n=0}^{N-1}\bar{\mathbf{x}}_{\alpha }\cdot
\nabla v_{n}\left( \mathbf{x}\right) =0.  \notag
\end{align}

Multiply both sides of \eqref{eq:3.4} by the function $\Psi _{m}(\alpha )$
for $0\leq m\leq N-1$ and then integrate the resulting equation with respect
to $\alpha $. We arrive at the following system of coupled quasilinear
elliptic equations: 
\begin{align}
& \Delta V\left( \mathbf{x}\right) +K\left( \nabla V\left( \mathbf{x}\right)
\right) =0\quad \text{for }\mathbf{x}\in \Omega ,  \label{eq:system} \\
& V\left( \mathbf{x}\right) =\psi _{0}\left( \mathbf{x}\right) \quad \text{%
for }\mathbf{x}\in \partial \Omega ,  \label{2} \\
& V_{z}\left( \mathbf{x}\right) =\psi _{1}\left( \mathbf{x}\right) \quad 
\text{for }\mathbf{x}\in \Gamma ,  \label{3} \\
& K\left( \nabla V\left( \mathbf{x}\right) \right) =S_{N}^{-1}f\left( \nabla
V\left( \mathbf{x}\right) \right) .  \label{202}
\end{align}%
Here $\psi _{0}\left( \mathbf{x}\right) $ and $\psi _{1}\left( \mathbf{x}%
\right) $ are known boundary data and we explain in \cref{sec:3.3} 3.2 how
to obtain them. Above, the unknown vector function $V(\mathbf{x})\in \mathbb{%
R}^{N}$ is given by 
\begin{equation*}
V\left( \mathbf{x}\right) =%
\begin{pmatrix}
v_{0}\left( \mathbf{x}\right) & v_{1}\left( \mathbf{x}\right) & \cdots & 
\cdots & v_{N-1}\left( \mathbf{x}\right)%
\end{pmatrix}%
^{T}.
\end{equation*}%
And $S_{N}=\left( s_{mn}\right) _{m,n=0}^{N-1}$ is the $N\times N$ matrix
that we have constructed above. The nonlinearity $f=\left( \left(
f_{m}\right) _{m=0}^{N-1}\right) ^{T}\in \mathbb{R}^{N}$ is quadratic with
respect to the first derivatives of components of $V\left( \mathbf{x}\right)
,$ 
\begin{align}
f_{m}\left( \nabla V\left( \mathbf{x}\right) \right) &
=2\sum_{n,l=0}^{N-1}\nabla v_{n}\left( \mathbf{x}\right) \cdot \nabla
v_{l}\left( \mathbf{x}\right) \int_{-a}^{a}\Psi _{m}\left( \alpha \right)
\Psi _{n}\left( \alpha \right) \Psi _{l}^{\prime }\left( \alpha \right)
d\alpha  \notag \\
& +2\sum_{n=0}^{N-1}\int_{-a}^{a}\Psi _{m}\left( \alpha \right) \Psi
_{n}^{\prime }\left( \alpha \right) \nabla v_{n}\left( \mathbf{x}\right)
\cdot \tilde{\mathbf{x}}_{\alpha }d\alpha  \label{203} \\
& +2\sum_{n=0}^{N-1}\int_{-a}^{a}\Psi _{m}\left( \alpha \right) \Psi
_{n}\left( \alpha \right) \bar{\mathbf{x}}_{\alpha }\cdot \nabla v_{n}\left( 
\mathbf{x}\right) d\alpha .  \notag
\end{align}%
It follows from (\ref{202}) and (\ref{203}) that the vector function $%
K\left( \nabla V\right) $ is quadratic with respect to components of $\nabla
V.$

The problem \eqref{eq:system}--(\ref{3}) is overdetermined since we have two
boundary conditions (\ref{2}), (\ref{3}) instead of just one. Also, this is
not a regular Cauchy problem for the system \eqref{eq:system} since the
Dirichlet data in (\ref{2}) are given at the entire boundary $\partial
\Omega .$ If solving problem \eqref{eq:system}--(\ref{3}), then we would
find the dielectric constant $c$ via backwards calculations. Therefore, we
focus below on the solution of problem \eqref{eq:system}--(\ref{3}).

\subsection{Boundary data (\protect\ref{2}), (\protect\ref{3})}

\label{sec:3.3} We now explain how to find the boundary data for the vector $%
V\left( \mathbf{x}\right) $ in (\ref{2}), (\ref{3}). It follows from %
\eqref{303} and (\ref{eq:truncated}) that the Dirichlet data at $\mathbf{x}%
\in \Gamma $ for $V\left( \mathbf{x}\right) $ are known. As it is known,
several data completion methods are heuristically applied in inverse
problems with incomplete data; see, e.g. \cite{Nguyen2018}. To complement
the lack of the boundary data information on $\partial \Omega \diagdown
\Gamma $, we use the data completion for \eqref{eq:aux}. More precisely, we
choose for each $\alpha $, 
\begin{equation}
\left. u\left( \mathbf{x},\alpha \right) \right\vert _{\partial \Omega }=%
\begin{cases}
F\left( \mathbf{x},\mathbf{x}_{\alpha }\right) , & \text{if }\mathbf{x}\in
\Gamma , \\ 
u_{0}\left( \mathbf{x},\alpha \right) , & \text{if }\mathbf{x}\in \partial
\Omega \diagdown \Gamma ,%
\end{cases}
\label{301}
\end{equation}%
where the $u_{0}\left( \mathbf{x},\alpha \right) $ is given in (\ref{1}) and
it is the solution of \eqref{eq:forward1}--\eqref{eq:forward2} for the case
of the uniform background. The choice \eqref{301} is fairly reasonable
because of the fact that $\partial \Omega \subset \left( \mathbb{R}%
^{3}\diagdown \Omega \right) $. Recall that the second condition (\ref{301})
is the final item of our approximate mathematical model (Remarks 3.1).

As to the data (\ref{3}), usually measurements are performed far from the
domain of interest, i.e. on the plane $\left\{ z=-R^{\prime }\right\} ,$
where $R^{\prime }>R.$ It is time consuming to solve a CIP in a large
domain. Besides, the data at the measurement plane are hard to use for an
inversion algorithm since they do not look \textquotedblleft nice". This can
be evidenced in our recent work for the experimental data; see, e.g., Figure
3a in \cite{JCP}. To \textquotedblleft move" the data closer to the target's
side, the so-called \textquotedblleft data propagation" procedure can be
applied to the measured data; see \cite{Nguyen2018} for a detailed
description of this procedure. By this procedure one obtains
\textquotedblleft propagated data", i.e. an approximation of the data at our
desired rectangle $\Gamma \subset \left\{ z=-R\right\} .$ Besides, the
propagated data look much better than the original data: e.g. compare
Figures 3a and 3b in \cite{JCP}. In addition, it is clear from the data
propagation procedure that one of its outcomes is an approximation of the $%
z- $derivative of the function $u\left( \mathbf{x},k\right) $ at $\Gamma .$
Thus, we assume that, in addition to the Dirichlet data at $\Gamma ,$ we
know the Neumann boundary data $u_{z}\left( \mathbf{x},\alpha \right)
=G\left( \mathbf{x},\alpha \right) $ for $\mathbf{x}\in \Gamma ,\mathbf{x}%
_{\alpha }\in L_{\text{src}}.$ Having the function $G\left( \mathbf{x}%
,\alpha \right) $ and using (\ref{eq:truncated}), one can easily find the
Neumann boundary condition $\psi _{1}\left( \mathbf{x}\right) $ at $\mathbf{x%
}\in \Gamma $ in (\ref{3}).

\subsection{Lipschitz stability of the boundary value problem (\protect\ref%
{eq:system})--(\protect\ref{3})}

\label{sec:3.4}

For any Banach space $B$ considered below and any integer $X>1$ we denote
the Banach space $B_{X}=\underbrace{B\times B\times ...\times B}_{X\text{
times}}$ with the norm%
\begin{equation*}
\left\Vert g\right\Vert _{B_{X}}=\left( \dsum\limits_{j=1}^{X}\left\Vert
g_{j}\right\Vert _{B}^{2}\right) ^{1/2}\quad \text{for all }g=\left(
g_{1},...,g_{K}\right) \in B_{X}.
\end{equation*}
\begin{verbatim}
 
\end{verbatim}

Let the number $r>R$ and the number $\lambda >0.$ In principle, many
functions can be used as CWFs for an elliptic operator. However, a rather
general one \cite{Klibanov2015} depends on two large parameters, which is
inconvenient for the numerical implementation. In our experience, better to
choose a rather simple Carleman Weight Function (CWF). Thus, we
define our CWF as%
\begin{equation}
\mu _{\lambda }\left( z\right) =\exp \left[ 2\lambda \left( z-r\right) ^{2}%
\right] ,z\in \left[ -R,R\right] .  \label{4}
\end{equation}%
We choose $r>R$ since one of conditions imposed on the CWF in any Carleman
estimate is that its gradient should not vanish in the closed domain.
Obviously, the function $\mu _{\lambda }\left( z\right) $ is decreasing for $%
z\in \left( -R,R\right) $ and 
\begin{equation}
\max_{\overline{\Omega }}\mu _{\lambda }\left( z\right) =\exp \left[
2\lambda \left( R+r\right) ^{2}\right] ,\text{ }\min_{\overline{\Omega }}\mu
_{\lambda }\left( z\right) =\exp \left[ 2\lambda \left( R-r\right) ^{2}%
\right] .  \label{40}
\end{equation}%
In other words, by (\ref{200}) and (\ref{201}) the CWF (\ref{4}) attains its
maximal value in $\overline{\Omega }$ on the part $\Gamma $ of the boundary
where measurements are conducted, and it attains its minimal value on the
opposite side.

Define the subspace $H_{0}^{2}\left( \Omega \right) $ of the space $%
H^{2}\left( \Omega \right) $ as: 
\begin{equation}
H_{0}^{2}\left( \Omega \right) :=\left\{ v\in H^{2}\left( \Omega \right)
:\left. v\right\vert _{\partial \Omega }=0,\left. \partial _{z}v\right\vert
_{\Gamma }=0\right\} .  \label{50}
\end{equation}%
Theorem 3.1 follows immediately from \cite[Theorem 4.1]{timedomain}.

\textbf{Theorem 3.1} (Carleman estimate). \emph{\ Let }$\mu _{\lambda
}\left( z\right) $\emph{\ be the function defined in (\ref{4}). Then there
exist constants }$\lambda _{0}=\lambda _{0}\left( \Omega ,r\right) \geq 1$%
\emph{\ and }$C=C\left( \Omega ,r\right) >0$\emph{\ depending only on the
domain }$\Omega $\emph{\ such that for every function }$u\in H_{0}^{2}\left(
\Omega \right) $\emph{\ and for all }$\lambda \geq \lambda _{0}$\emph{\ the
following Carleman estimate holds:}%
\begin{equation}
\dint\limits_{\Omega }\left\vert \Delta u\right\vert ^{2}\mu _{\lambda
}\left( z\right) d\mathbf{x}\geq \frac{C}{\lambda }\dsum\limits_{i,j=1}^{3}%
\dint\limits_{\Omega }\left\vert u_{x_{i}x_{j}}\right\vert ^{2}\mu _{\lambda
}\left( z\right) d\mathbf{x}+C\lambda \dint\limits_{\Omega }\left[
\left\vert \nabla u\right\vert ^{2}+\lambda ^{2}\left\vert u\right\vert ^{2}%
\right] \mu _{\lambda }\left( z\right) d\mathbf{x.}  \label{6}
\end{equation}

Suppose that there exist two vector functions $V^{\left( 1\right) }\left( 
\mathbf{x}\right) $ and $V^{\left( 2\right) }\left( \mathbf{x}\right) $
satisfying equation (\ref{eq:system}) with boundary conditions as in (\ref{2}%
), (\ref{3}), 
\begin{align}
& V^{\left( 1\right) }\left( \mathbf{x}\right) =\psi _{0}^{\left( 1\right)
}\left( \mathbf{x}\right) ,V^{\left( 2\right) }\left( \mathbf{x}\right)
=\psi _{0}^{\left( 2\right) }\left( \mathbf{x}\right) \quad \text{for }%
\mathbf{x}\in \partial \Omega ,  \label{7} \\
& V_{z}^{\left( 1\right) }\left( \mathbf{x}\right) =\psi _{1}^{\left(
1\right) }\left( \mathbf{x}\right) ,V_{z}^{\left( 2\right) }\left( \mathbf{x}%
\right) =\psi _{1}^{\left( 2\right) }\left( \mathbf{x}\right) \quad \text{%
for }\mathbf{x}\in \Gamma .  \label{8}
\end{align}%
Suppose that there exist two vector functions $F_{1}$,$F_{2}\in
H_{N}^{3}\left( \Omega \right) $ satisfying boundary conditions (\ref{7}), (%
\ref{8}), i.e.%
\begin{align}
& F_{1}\left( \mathbf{x}\right) =\psi _{0}^{\left( 1\right) }\left( \mathbf{x%
}\right) ,F_{2}\left( \mathbf{x}\right) =\psi _{0}^{\left( 2\right) }\left( 
\mathbf{x}\right) ,\quad \text{for }\mathbf{x}\in \partial \Omega ,
\label{9} \\
& \partial _{z}F_{1}\left( \mathbf{x}\right) =\psi _{1}^{\left( 1\right)
}\left( \mathbf{x}\right) ,\partial _{z}F_{2}\left( \mathbf{x}\right) =\psi
_{1}^{\left( 2\right) }\left( \mathbf{x}\right) ,\quad \text{for }\mathbf{x}%
\in \Gamma .  \label{10}
\end{align}%
Let $M>0$ be a number. We assume that%
\begin{equation}
V^{\left( 1\right) },V^{\left( 2\right) },F_{1},F_{2}\in G\left( M\right)
=\left\{ W\in H_{N}^{3}\left( \Omega \right) :\left\Vert W\right\Vert
_{H_{N}^{3}\left( \Omega \right) }<M\right\} .  \label{11}
\end{equation}%
Note also that by the embedding theorem 
\begin{equation}
G\left( M\right) \subset C_{N}^{1}\left( \overline{\Omega }\right) \text{
and }\left\Vert W\right\Vert _{C_{N}^{1}\left( \overline{\Omega }\right)
}\leq C_{1}\quad \text{for all }W\in G\left( M\right) .  \label{12}
\end{equation}%
Here and below $C_{1}=C_{1}\left( \Omega ,N,M\right) >0$ denotes different
constants depending only on listed parameters.

\textbf{Theorem 3.2} (Lipschitz stability estimate). \emph{Let }$V^{\left(
1\right) }\left( \mathbf{x}\right) $\emph{\ and }$V^{\left( 2\right) }\left( 
\mathbf{x}\right) $\emph{\ be two solutions of equation (\ref{eq:system})
with boundary conditions (\ref{7}), (\ref{8}). Suppose that there exist two
vector functions }$F_{1} $\emph{,}$F_{2} \in H_{N}^{3}\left( \Omega \right) $%
\emph{\ satisfying (\ref{9}), (\ref{10}). Also, let (\ref{11}) holds. Then
the following Lipschitz stability estimate is valid}%
\begin{equation}
\left\Vert V^{\left( 1\right) }-V^{\left( 2\right) }\right\Vert
_{H_{N}^{2}\left( \Omega \right) }\leq C_{1}\left\Vert
F_{1}-F_{2}\right\Vert _{H_{N}^{2}\left( \Omega \right) }.  \label{120}
\end{equation}

\textbf{Proof}. Denote 
\begin{align}
& Q_{1}\left( \mathbf{x}\right) =V^{\left( 1\right) }\left( \mathbf{x}%
\right) -F_{1}\left( \mathbf{x}\right) ,\quad Q_{2}\left( \mathbf{x}\right)
=V^{\left( 2\right) }\left( \mathbf{x}\right) -F_{2}\left( \mathbf{x}\right)
,  \label{13} \\
& \widetilde{Q}\left( \mathbf{x}\right) =Q_{1}\left( \mathbf{x}\right)
-Q_{2}\left( \mathbf{x}\right) ,\quad \widetilde{F}\left( \mathbf{x}\right)
=F_{1}\left( \mathbf{x}\right) -F_{2}\left( \mathbf{x}\right) .  \label{14}
\end{align}%
Then (\ref{203}) and (\ref{7})-(\ref{14}) imply that%
\begin{align}
& \Delta \widetilde{Q}\left( \mathbf{x}\right) =T_{1}\left( \mathbf{x}%
\right) \cdot \nabla \widetilde{Q}\left( \mathbf{x}\right) +T_{2}\left( 
\mathbf{x}\right) \cdot \nabla \widetilde{F}\left( \mathbf{x}\right) -\Delta 
\widetilde{F}\left( \mathbf{x}\right) ,  \label{15} \\
& \widetilde{Q}\mid _{\partial \Omega }=0,\quad \widetilde{Q}_{z}\mid
_{\Gamma }=0,  \label{16} \\
& T_{1},T_{2}\in C_{N}^{3}\left( \overline{\Omega }\right) ,\left\Vert
T_{1}\right\Vert _{C_{N}^{3}\left( \overline{\Omega }\right) },\left\Vert
T_{2}\right\Vert _{C_{N}^{3}\left( \overline{\Omega }\right) }\leq C_{1}.
\label{17}
\end{align}%
Square absolute values of both sides of equation (\ref{15}). Next, multiply
the resulting equation by the CWF (\ref{4}) and integrate over the domain $%
\Omega .$ Using (\ref{17}), we obtain 
\begin{equation}
\dint\limits_{\Omega }\left\vert \Delta \widetilde{Q}\right\vert ^{2}\mu
_{\lambda }\left( z\right) d\mathbf{x}\leq C_{1}\dint\limits_{\Omega
}\left\vert \nabla \widetilde{Q}\right\vert ^{2}\mu _{\lambda }\left(
z\right) d\mathbf{x}+C_{1}\dint\limits_{\Omega }\left( \left\vert \Delta 
\widetilde{F}\right\vert ^{2}+\left\vert \nabla \widetilde{F}\right\vert
^{2}\right) \mu _{\lambda }\left( z\right) d\mathbf{x.}  \label{18}
\end{equation}%
Taking into account (\ref{50}) and (\ref{16}) and also applying (\ref{6}) to
(\ref{18}), we obtain for all $\lambda \geq \lambda _{0}>1$ 
\begin{align}
& C_{1}\dint\limits_{\Omega }\left( \left\vert \Delta \widetilde{F}%
\right\vert ^{2}+\left\vert \nabla \widetilde{F}\right\vert ^{2}\right) \mu
_{\lambda }\left( z\right) d\mathbf{x+}C_{1}\dint\limits_{\Omega }\left\vert
\nabla \widetilde{Q}\right\vert ^{2}\mu _{\lambda }\left( z\right) d\mathbf{x%
}  \label{19} \\
& \geq \frac{1}{\lambda }\dsum\limits_{i,j=1}^{3}\dint\limits_{\Omega
}\left\vert \widetilde{Q}_{x_{i}x_{j}}\right\vert ^{2}\mu _{\lambda }\left(
z\right) d\mathbf{x}+\lambda \dint\limits_{\Omega }\left[ \left\vert \nabla 
\widetilde{Q}\right\vert ^{2}+\lambda ^{2}\left\vert \widetilde{Q}%
\right\vert ^{2}\right] \mu _{\lambda }\left( z\right) d\mathbf{x}.  \notag
\end{align}%
Choose a number $\lambda _{1}\geq \lambda _{0}$ such that $\lambda _{1}\geq
2C_{1}.$ Then (\ref{19}) implies that%
\begin{align*}
& C_{1}\dint\limits_{\Omega }\left( \left\vert \Delta \widetilde{F}%
\right\vert ^{2}+\left\vert \nabla \widetilde{F}\right\vert ^{2}\right) \mu
_{\lambda _{1}}\left( z\right) d\mathbf{x} \\
& \geq \frac{1}{\lambda _{1}}\dsum\limits_{i,j=1}^{3}\dint\limits_{\Omega
}\left\vert \widetilde{Q}_{x_{i}x_{j}}\right\vert ^{2}\mu _{\lambda
_{1}}\left( z\right) d\mathbf{x}+\frac{\lambda _{1}}{2}\dint\limits_{\Omega }%
\left[ \left\vert \nabla \widetilde{Q}\right\vert ^{2}+\left\vert \widetilde{%
Q}\right\vert ^{2}\right] \mu _{\lambda _{1}}\left( z\right) d\mathbf{x}.
\end{align*}%
This inequality and (\ref{40}) lead to:%
\begin{align*}
& C_{1}\exp \left( 4Rr\lambda _{1}\right) \dint\limits_{\Omega }\left(
\left\vert \Delta \widetilde{F}\right\vert ^{2}+\left\vert \nabla \widetilde{%
F}\right\vert ^{2}\right) d\mathbf{x} \\
& \geq \frac{1}{\lambda _{1}}\dsum\limits_{i,j=1}^{3}\dint\limits_{\Omega
}\left\vert \widetilde{Q}_{x_{i}x_{j}}\right\vert ^{2}d\mathbf{x}+\frac{%
\lambda _{1}}{2}\dint\limits_{\Omega }\left[ \left\vert \nabla \widetilde{Q}%
\right\vert ^{2}+\left\vert \widetilde{Q}\right\vert ^{2}\right] d\mathbf{x}.
\end{align*}%
Hence, with a new constant $C_{1}$ we have%
\begin{equation}
\left\Vert \widetilde{Q}\right\Vert _{H_{N}^{2}\left( \Omega \right) }\leq
C_{1}\left\Vert \widetilde{F}\right\Vert _{H_{N}^{2}\left( \Omega \right) }.
\label{20}
\end{equation}%
Next, by (\ref{13}), (\ref{14}) and triangle inequality%
\begin{align*}
\left\Vert \widetilde{Q}\right\Vert _{H_{N}^{2}\left( \Omega \right) }&
=\left\Vert \left( V^{\left( 1\right) }-F_{1}\right) -\left( V^{\left(
2\right) }-F_{2}\right) \right\Vert _{H_{N}^{2}\left( \Omega \right) } \\
& \geq \left\Vert V^{\left( 1\right) }-V^{\left( 2\right) }\right\Vert
_{H_{N}^{2}\left( \Omega \right) }-\left\Vert F_{1}-F_{2}\right\Vert
_{H_{N}^{2}\left( \Omega \right) }.
\end{align*}%
Combining this with (\ref{20}), we obtain the target estimate (\ref{120}) of
this theorem. $\square $

\section{Weighted Tikhonov-like Functional}

\label{sec:4}

For the convenience of the presentation, each $N-$D complex valued vector
function $W=\func{Re}W+i\func{Im}W$ is considered below as the $2N-$D vector
function with real valued components $\left( \func{Re}W,\func{Im}W\right)
:=\left( W_{1},W_{2}\right) :=W\in \mathbb{R}^{2N}$. All results and proofs
below are for these $2N-$D vector functions. For any number $s\in \mathbb{C}$%
, its complex conjugate is denoted as $\overline{s}$.

We find an approximate solution of the problem \eqref{eq:system}-(\ref{202})
via the minimization of an appropriate weighted Tikhonov-like functional
with the CWF (\ref{4}) involved in it. Due to \eqref{eq:system}, denote 
\begin{equation}
L\left( V\right) \left( \mathbf{x}\right) =\Delta V\left( \mathbf{x}\right)
+K\left( \nabla V\left( \mathbf{x}\right) \right) .  \label{4.1}
\end{equation}%
Let $\gamma \in (0,1)$ be the regularization parameter. We now consider the
following weighted Tikhonov-like functional $J_{\lambda ,\gamma
}:H_{2N}^{3}(\Omega )\rightarrow \mathbb{R}_{+}$, 
\begin{equation}
J_{\lambda ,\gamma }\left( V\right) =\exp \left[ -2\lambda \left( R+r\right)
^{2}\right] \dint\limits_{\Omega }\left\vert L\left( V\right) \right\vert
^{2}\mu _{\lambda }\left( z\right) d\mathbf{x}+\gamma \left\Vert
V\right\Vert _{H_{N}^{3}\left( \Omega \right) }^{2}.  \label{4.2}
\end{equation}%
Here $\exp \left[ -2\lambda \left( R+r\right) ^{2}\right] $ is the balancing
multiplier: to balance first and second terms in the right hand side of (\ref%
{4.2}), see (\ref{40}). We use the $H_{N}^{3}\left( \Omega \right) -$norm in
the regularization term here since $H_{N}^{3}\left( \Omega \right) \subset
C_{N}^{1}\left( \overline{\Omega }\right) $ and an obvious analog of (\ref%
{12}) holds.

Assuming for a moment that the nonlinear term $K\left( \nabla V\left( 
\mathbf{x}\right) \right) $ is absent in (\ref{4.1}), we remark that since
the Laplace operator is linear, then one can also find an approximate
solution of the problem \eqref{eq:system}--(\ref{202}) by the regular
quasi-reversibility method with $\lambda =0$ in (\ref{4.2}) (see, e.g. \cite%
{Klibanov2015}). However, if $K\left( \nabla V\left( \mathbf{x}\right)
\right) \neq 0,$ then the presence of the CWF serves three purposes: first,
it controls this nonlinear term; second, it \textquotedblleft
maximizes\textquotedblright\ the influence of the important boundary data at 
$z=-R$; and third, it \textquotedblleft convexifies" the cost functional
globally. These are the underlying reasons of the convexification idea.
Below $\left( \cdot ,\cdot \right) $ is the scalar product in the space $%
H_{2N}^{3}(\Omega )$.\textbf{\ }Let $M>0$ be an arbitrary number. We define
the set $B\left( M\right) \subset H_{2N}^{3}(\Omega )$ as%
\begin{equation}
B\left( M\right) =\left\{ 
\begin{array}{c}
V\in H_{2N}^{3}(\Omega ):\left\Vert V\right\Vert _{H_{2N}^{3}(\Omega
)}<M,V\mid _{\partial \Omega }=\psi _{0},V_{z}\mathbf{\mid }_{\Gamma }=\psi
_{1}%
\end{array}%
\right\} .  \label{4.03}
\end{equation}%
By (\ref{12}), we know that%
\begin{equation}
B\left( M\right) \subset C_{2N}^{1}\left( \overline{\Omega }\right) \text{
and }\left\Vert V\right\Vert _{C_{2N}^{1}\left( \overline{\Omega }\right)
}\leq C_{1}\quad \text{ for all }V\in B\left( M\right) .  \label{4.3}
\end{equation}

\textbf{Minimization problem (MP). }\emph{Minimize the cost functional }$%
J_{\lambda ,\gamma }(V)$\emph{\ on the set }$\overline{B\left( M\right) }.$

\section{Analysis of the Functional $J_{\protect\lambda ,\protect\gamma %
}\left( V\right) $}

\label{sec:5}

\subsection{Strict \ convexity on $\overline{B\left( M\right) }$}

\label{sec:5.1}

Theorem 5.1 is the central analytical result of this work. Note that in the
proof of this theorem we do not \textquotedblleft subtract" boundary
conditions from the vector function $V$, which means that we do not arrange
zero boundary conditions for the difference. Hence, we do not require here
that our boundary conditions should be extended in the entire domain $\Omega
.$ This is a new element compared with our proofs in the previous works on
the convexification \cite{Bak,convexper,convIPnew,EIT,timedomain}.
Thus, the proof of Theorem 5.1 is significantly different from those of
these works. Still, we use that subtraction in Theorems 5.4 and 6.1.

\textbf{Theorem 5.1}. \emph{The functional }$J_{\lambda ,\gamma }\left(
V\right) $\emph{\ has its Frech\'{e}t derivative }$J_{\lambda ,\gamma
}^{\prime }\left( V\right) $\emph{\ at any point }$V\in \overline{B\left(
M\right) }.$\emph{\ Let }$\lambda _{0}>1$\emph{\ be the number of Theorem
3.1. There exists a sufficiently large number }$\lambda _{2}=\lambda
_{2}\left( M,N,r,\Omega \right) \geq \lambda _{0}$\emph{\ such that the
functional }$J_{\lambda ,\gamma }\left( V\right) $\emph{\ is strictly convex
on }$\overline{B\left( M\right) }$\emph{\ for all }$\lambda \geq \lambda
_{2}.$\emph{\ More precisely, for all $\lambda \geq \lambda _{2}$ the
following inequality holds}%
\begin{align}
& J_{\lambda ,\gamma }\left( V^{\left( 2\right) }\right) -J_{\lambda ,\gamma
}\left( V^{\left( 1\right) }\right) -J_{\lambda ,\gamma }^{\prime }\left(
V^{\left( 1\right) }\right) \left( V^{\left( 2\right) }-V^{\left( 1\right)
}\right)  \label{5.1} \\
\geq C_{1}& \left\Vert V^{\left( 2\right) }-V^{\left( 1\right) }\right\Vert
_{H_{2N}^{2}\left( \Omega \right) }^{2}+\gamma \left\Vert V^{\left( 2\right)
}-V^{\left( 1\right) }\right\Vert _{H_{2N}^{3}\left( \Omega \right)
}^{2}\;\text{ for all } V^{\left( 1\right) },V^{\left( 2\right) }\in \overline{%
B\left( M\right) }.  \notag
\end{align}

\textbf{Proof}. Let $V^{\left( 1\right) },V^{\left( 2\right) }\in \overline{%
B\left( M\right) }$ be two arbitrary points. Define $H_{0,2N}^{3}\left(
\Omega \right) =H_{2N}^{3}\left( \Omega \right) \cap H_{0,2N}^{2}\left(
\Omega \right) ;$ see (\ref{50}). Denote $h=\left( h_{1},h_{2}\right)
=V^{\left( 2\right) }-V^{\left( 1\right) }.$ Then%
\begin{equation}
h\in \overline{B\left( 2M\right) }\;\text{ and }\;h\in H_{0,2N}^{3}\left(
\Omega \right) .  \label{5.2}
\end{equation}
Obviously $\left\vert L\left( V^{\left( 2\right) }\right)
\right\vert ^{2}=\left\vert L\left( V^{\left( 1\right) }+h\right)
\right\vert ^{2}.$ Observe that it follows from (\ref{eq:system}), (\ref{202}) and (\ref{203}) that the vector function $K\left( \nabla
V\right) $ is the sum of linear and quadratic parts with respect to
the gradients $\nabla v_{n}\left( \mathbf{x}\right) $ of the
components $v_{n}\left( \mathbf{x}\right) $ of the vector function $V.$ Using this as well as (\ref{4.1}), we obtain%
\begin{equation}
L\left( V^{\left( 1\right) }+h\right) =L\left( V^{\left( 1\right) }\right)
+\Delta h+K_{1}\left( \mathbf{x}\right) \nabla h+K_{2}\left( \mathbf{x}%
,\nabla h\right) .  \label{5.3}
\end{equation}%
Here, the vector functions $K_{1},K_{2}$ are continuous with respect to $%
\mathbf{x}$ in $\overline{\Omega }.$ 
Also, $K_{1}\left( \mathbf{x}\right) $ is independent on $h$. As to the vector function 
$K_{2}\left( \mathbf{x},\nabla h\right) ,$ it is quadratic with
respect to the gradients $\nabla h_{n}\left( \mathbf{x}\right) $ 
of the components of the vector function $h.$ The latter, (\ref{4.03}) and (\ref{4.3}) imply that
\begin{equation}
\left\vert K_{2}\left( \mathbf{x},\nabla h\right) \right\vert \leq
C_{1}\left\vert \nabla h\right\vert ^{2}\quad \text{for all }\mathbf{x}\in 
\overline{\Omega }.  \label{5.4}
\end{equation}

Squaring absolute values of both sides of (\ref{5.3}), we
obtain %
\begin{align}
\left\vert L\left( V^{\left( 1\right) }+h\right) \right\vert
^{2}&=\left\vert L\left( V^{\left( 1\right) }\right) \right\vert ^{2}+2\func{%
Re}\left\{ \overline{L\left( V^{\left( 1\right) }\right) }\left[ \Delta
h+K_{1}\left( \mathbf{x}\right) \nabla h+K_{2}\left( \mathbf{x},\nabla
h\right) \right] \right\}  \notag \\
& +\left\vert \Delta h+K_{1}\left( \mathbf{x}\right) \nabla h+K_{2}\left( 
\mathbf{x},\nabla h\right) \right\vert ^{2}.  \label{5.40}
\end{align}%
In (\ref{5.40}), we single out the linear, with respect to
$h$, term as well as the term $\left\vert \Delta h\right\vert ^{2}$. We
obtain%
\begin{align}
& \left\vert L\left( V^{\left( 1\right) }+h\right) \right\vert
^{2}-\left\vert L\left( V^{\left( 1\right) }\right) \right\vert ^{2}
\label{5.50} \\
& =2\func{Re}\left\{ \overline{L\left( V^{\left( 1\right) }\right) }\left[
\Delta h+K_{1}\left( \mathbf{x}\right) \nabla h\right] \right\} +\left\vert
\Delta h\right\vert ^{2} +2\func{Re}\left[ \overline{L\left( V^{\left(
1\right) }\right) }K_{2}\left( \mathbf{x},\nabla h\right) \right]  \notag \\
& +2\func{Re}\left\{ \overline{\Delta h}\left[ K_{1}\left( \mathbf{x}\right)
\nabla h+K_{2}\left( \mathbf{x},\nabla h\right) \right] \right\} +\left\vert
K_{1}\left( \mathbf{x}\right) \nabla h+K_{2}\left( \mathbf{x},\nabla
h\right) \right\vert ^{2}.  \notag
\end{align}%
In (\ref{5.50}), the term $2\func{Re}\left\{ \overline{L\left( V^{\left( 1\right) }\right) }\left[
	\Delta h+K_{1}\left( \mathbf{x}\right) \nabla h\right] \right\} $ is linear with respect to $h$. Thus, we obtain
\begin{align}  \label{5.5}
& J_{\lambda ,\gamma }\left( V^{\left( 1\right) }+h\right) -J_{\lambda
,\gamma }\left( V^{\left( 1\right) }\right) =Lin\left( h\right) +\gamma
\left\Vert h\right\Vert _{H_{2N}^{3}\left( \Omega \right) }^{2} \\
& +e^{-2\lambda \left( R+r\right) ^{2}}\dint\limits_{\Omega }\left\{
\left\vert \Delta h\right\vert ^{2}+2\func{Re}\left[ \overline{\Delta h}%
\cdot \left( K_{1}\left( \mathbf{x}\right) \nabla h+K_{2}\left( \mathbf{x}%
,\nabla h\right) \right) \right] \right\} \mu _{\lambda }\left( z\right) d%
\mathbf{x}  \notag \\
& +e^{-2\lambda \left( R+r\right) ^{2}}\dint\limits_{\Omega }\left[ 2\func{Re%
}\left[ \overline{L\left( V^{\left( 1\right) }\right) }K_{2}\left( \mathbf{x}%
,\nabla h\right) \right] +\left\vert K_{1}\left( \mathbf{x}\right) \nabla
h+K_{2}\left( \mathbf{x},\nabla h\right) \right\vert ^{2}\right] \mu
_{\lambda }\left( z\right) d\mathbf{x,}  \notag
\end{align}%
where the functional $Lin\left( h\right) :H_{0,2N}^{3}\rightarrow \mathbb{R}$
is linear with respect to $h=\left( h_{1},h_{2}\right) $. 
It follows
from (\ref{5.50}) that it is generated by the term $2\func{Re}\left\{ 
\overline{L\left( V^{\left( 1\right) }\right) }\left[ \Delta h+K_{1}\left( 
\mathbf{x}\right) \nabla h\right] \right\} $, 
\begin{align}  \label{5.6}
& Lin\left( h\right) = 2\gamma \left( V^{\left( 1\right) },h\right) \\
& +2e^{-2\lambda \left( R+r\right) ^{2}}\dint\limits_{\Omega }\func{Re}\left[
\Delta h+\left( K_{1}\left( \mathbf{x}\right) \nabla h\right) \left( 
\overline{\Delta V^{\left( 1\right) }+K\left( \nabla V^{\left( 1\right)
}\right) }\right) \right] \mu _{\lambda }\left( z\right) d\mathbf{x} . 
\notag
\end{align}%
Besides, it follows from (\ref{5.5}) and (\ref{5.6}) that 
\begin{equation*}
\lim_{\left\Vert h\right\Vert _{H_{2N}^{3}\left( \Omega \right) }\rightarrow
0^{+}}\left\{ \frac{1}{\left\Vert h\right\Vert _{H_{2N}^{3}\left( \Omega
\right) }}\left[ J_{\lambda ,\gamma }\left( V^{\left( 1\right) }+h\right)
-J_{\lambda ,\gamma }\left( V^{\left( 1\right) }\right) -Lin\left( h\right) %
\right] \right\} =0.
\end{equation*}%
Hence, the functional $Lin\left( h\right) $ is the Frech\'{e}t derivative of
the functional $J_{\lambda ,\gamma }$ at the point $V^{\left( 1\right) }\in 
\overline{B\left( M\right) }.$ By the Riesz theorem, there exists a unique
point $J_{\lambda ,\gamma }^{\prime }\left( V^{\left( 1\right) }\right) $
such that 
\begin{align}\label{555}
J_{\lambda ,\gamma }^{\prime }\left( V^{\left( 1\right) }\right) \in
H_{0,2N}^{3}\left( \Omega \right)  \text{ and } Lin\left( h\right) =\left(
J_{\lambda ,\gamma }^{\prime }\left( V^{\left( 1\right) }\right) ,h\right)
\quad \text{for all }h\in H_{0,2N}^{3}\left( \Omega \right) .
\end{align}
Thus, we can rewrite (\ref{5.5}) as%
\begin{align}  \label{5.7}
& J_{\lambda ,\gamma }\left( V^{\left( 1\right) }+h\right) -J_{\lambda
,\gamma }\left( V^{\left( 1\right) }\right) -\left( J_{\lambda ,\gamma
}^{\prime }\left( V^{\left( 1\right) }\right) ,h\right) =\gamma \left\Vert
h\right\Vert _{H_{2N}^{3}\left( \Omega \right) }^{2} \\
& +e^{-2\lambda \left( R+r\right) ^{2}}\dint\limits_{\Omega }\left\{
\left\vert \Delta h\right\vert ^{2}+2\func{Re}\left[ \overline{\Delta h}%
\cdot \left( K_{1}\left( \mathbf{x}\right) \nabla h+K_{2}\left( \mathbf{x}%
,\nabla h\right) \right) \right] \right\} \mu _{\lambda }\left( z\right) d%
\mathbf{x}  \notag \\
& +e^{-2\lambda \left( R+r\right) ^{2}}\dint\limits_{\Omega }\left[ 2\func{Re%
}\left[ \overline{L\left( V^{\left( 1\right) }\right) }K_{2}\left( \mathbf{x}%
,\nabla h\right) \right] +\left\vert K_{1}\left( \mathbf{x}\right) \nabla
h+K_{2}\left( \mathbf{x},\nabla h\right) \right\vert ^{2}\right] \mu
_{\lambda }\left( z\right) d\mathbf{x.}  \notag
\end{align}%
We now estimate from the below the term in the second line of (\ref{5.7}).
By the Cauchy--Schwarz inequality, (\ref{4.3}) and (\ref{5.4}) we find that 
\begin{equation*}
2\left\vert \overline{\Delta h}\cdot \left( K_{1}\left( \mathbf{x}\right)
\nabla h+K_{2}\left( \mathbf{x},\nabla h\right) \right) \right\vert \leq 
\frac{1}{2}\left\vert \Delta h\right\vert ^{2}+C_{1}\left\vert \nabla
h\right\vert ^{2}.
\end{equation*}%
Therefore, 
\begin{align}
& \dint\limits_{\Omega }\left\{ \left\vert \Delta h\right\vert ^{2}+2\func{Re%
}\left[ \overline{\Delta h}\cdot \left( K_{1}\left( \mathbf{x}\right) \nabla
h+K_{2}\left( \mathbf{x},\nabla h\right) \right) \right] \right\} \mu
_{\lambda }\left( z\right) d\mathbf{x}  \notag \\
& \geq \dint\limits_{\Omega }\left\vert \Delta h\right\vert ^{2}\mu
_{\lambda }\left( z\right) d\mathbf{x}-\frac{1}{2}\dint\limits_{\Omega
}\left\vert \Delta h\right\vert ^{2}\mu _{\lambda }\left( z\right) d\mathbf{x%
}-C_{1}\dint\limits_{\Omega }\left\vert \nabla h\right\vert ^{2}\mu
_{\lambda }\left( z\right) d\mathbf{x}  \label{5.8} \\
& =\frac{1}{2}\dint\limits_{\Omega }\left\vert \Delta h\right\vert ^{2}\mu
_{\lambda }\left( z\right) d\mathbf{x}-C_{1}\dint\limits_{\Omega }\left\vert
\nabla h\right\vert ^{2}\mu _{\lambda }\left( z\right) d\mathbf{x.}  \notag
\end{align}%
Next, using (\ref{5.4}), we estimate from the below the term in the third
line of (\ref{5.7}),%
\begin{equation*}
e^{-2\lambda \left( R+r\right) ^{2}}\dint\limits_{\Omega }\left[ 2\func{Re}%
\left[ \overline{L\left( V^{\left( 1\right) }\right) }K_{2}\left( \mathbf{x}%
,\nabla h\right) \right] +\left\vert K_{1}\left( \mathbf{x}\right) \nabla
h+K_{2}\left( \mathbf{x},\nabla h\right) \right\vert ^{2}\right] \mu
_{\lambda }\left( z\right) d\mathbf{x}
\end{equation*}%
\begin{equation}
\geq -C_{1}e^{-2\lambda \left( R+r\right) ^{2}}\dint\limits_{\Omega
}\left\vert \nabla h\right\vert ^{2}\mu _{\lambda }\left( z\right) d\mathbf{%
x.}  \label{5.80}
\end{equation}%
Thus, (\ref{5.7})--(\ref{5.80}) imply

\begin{align}
& J_{\lambda ,\gamma }\left( V^{\left( 1\right) }+h\right) -J_{\lambda
,\gamma }\left( V^{\left( 1\right) }\right) -\left( J_{\lambda ,\gamma
}^{\prime }\left( V^{\left( 1\right) }\right) ,h\right)  \label{5.9} \\
& \geq \frac{e^{-2\lambda \left( R+r\right) ^{2}}}{2}\left[
\dint\limits_{\Omega }\left\vert \Delta h\right\vert ^{2}\mu _{\lambda
}\left( z\right) d\mathbf{x}-C_{1}\dint\limits_{\Omega }\left\vert \nabla
h\right\vert ^{2}\mu _{\lambda }\left( z\right) d\mathbf{x}\right] +\gamma
\left\Vert h\right\Vert _{H_{2N}^{3}\left( \Omega \right) }^{2}.  \notag
\end{align}%
Now we apply the Carleman estimate (\ref{6}) to the second line of (\ref{5.9}%
). This use is possible due to (\ref{5.2}). For brevity, we do not count the
multiplier $\exp \left[ -2\lambda \left( R+r\right) ^{2}\right] $ for a
while. With a constant $\widetilde{C}=\widetilde{C}\left( \Omega ,r,N\right)
>0$ and a number $\widetilde{\lambda }_{0}=\widetilde{\lambda }_{0}\left(
\Omega ,r,N\right) \geq \lambda _{0}>1$ depending only on listed parameters,
we obtain for all $\lambda \geq \widetilde{\lambda }_{0}$ 
\begin{align}
& \frac{1}{2}\dint\limits_{\Omega }\left\vert \Delta h\right\vert ^{2}\mu
_{\lambda }\left( z\right) d\mathbf{x}-C_{1}\dint\limits_{\Omega }\left\vert
\nabla h\right\vert ^{2}\mu _{\lambda }\left( z\right) d\mathbf{x}\geq \frac{%
\widetilde{C}}{\lambda }\dsum\limits_{i,j=1}^{3}\dint\limits_{\Omega
}\left\vert h_{x_{i}x_{j}}\right\vert ^{2}\mu _{\lambda }\left( z\right) d%
\mathbf{x}  \label{5.10} \\
& +\widetilde{C}\lambda \dint\limits_{\Omega }\left[ \left\vert \nabla
h\right\vert ^{2}+\lambda ^{2}\left\vert h\right\vert ^{2}\right] \mu
_{\lambda }\left( z\right) d\mathbf{x}-C_{1}\dint\limits_{\Omega }\left\vert
\nabla h\right\vert ^{2}\mu _{\lambda }\left( z\right) d\mathbf{x}.  \notag
\end{align}%
Choose the number $\lambda _{2}=\lambda _{2}\left( M,\Omega ,r,N\right) \geq 
\widetilde{\lambda }_{0}>1$ depending only on listed parameters such that $%
\widetilde{C}\lambda _{2}>2C_{1}.$ Then we obtain from (\ref{5.10}) 
\begin{align}
& \frac{1}{2}\dint\limits_{\Omega }\left\vert \Delta h\right\vert ^{2}\mu
_{\lambda _{2}}\left( z\right) d\mathbf{x}-C_{1}\dint\limits_{\Omega
}\left\vert \nabla h\right\vert ^{2}\mu _{\lambda _{2}}\left( z\right) d%
\mathbf{x}  \notag \\
& \geq \frac{\widetilde{C}}{\lambda _{2}}\dsum\limits_{i,j=1}^{3}\dint%
\limits_{\Omega }\left\vert h_{x_{i}x_{j}}\right\vert ^{2}\mu _{\lambda
_{2}}\left( z\right) d\mathbf{x}+\frac{1}{2}\widetilde{C}\lambda
_{2}\dint\limits_{\Omega }\left[ \left\vert \nabla h\right\vert ^{2}+\lambda
_{2}^{2}\left\vert h\right\vert ^{2}\right] \mu _{\lambda _{2}}\left(
z\right) d\mathbf{x}  \label{5.11} \\
& \geq C_{1}e^{2\lambda _{2}\left( R-r\right) ^{2}}\left\Vert h\right\Vert
_{H_{2N}^{2}}^{2}.  \notag
\end{align}%
Hence, combining (\ref{5.9})--(\ref{5.11}) we arrive at 
\begin{equation*}
J_{\lambda ,\gamma }\left( V^{\left( 1\right) }+h\right) -J_{\lambda ,\gamma
}\left( V^{\left( 1\right) }\right) -\left( J_{\lambda ,\gamma }^{\prime
}\left( V^{\left( 1\right) }\right) ,h\right) \geq C_{1}\left\Vert
h\right\Vert _{H_{2N}^{2}}^{2}+\gamma \left\Vert h\right\Vert
_{H_{2N}^{3}\left( \Omega \right) }^{2},
\end{equation*}%
which is equivalent to our target estimate (\ref{5.1}). \ \ \ \ $\square $

\subsection{The minimizer of $J_{\protect\lambda ,\protect\gamma }\left(
V\right) $ on $\overline{B\left( M\right) }$}

\label{sec:5.2}

In Theorem 5.2 below, we state the Lipschitz continuity of the Frech\'{e}t
derivative $J_{\lambda ,\gamma }^{\prime }\left( V\right) $ on $\overline{%
B\left( M\right) }.$ We omit the proof of this theorem because it is very
similar to the proof of Theorem 3.1 in \cite{Bak}.

\textbf{Theorem 5.2}. \emph{For any }$\lambda >0$\emph{\ the Frech\'{e}t
derivative }$J_{\lambda ,\gamma }^{\prime }\left( V\right) $\emph{\ of the
functional }$J_{\lambda ,\gamma }\left( V\right) $\emph{\ is Lipschitz
continuous on the set }$\overline{B\left( M\right) }.$\emph{\ In other
words, there exists a number }$D=D\left( \Omega ,r,N,M,\lambda ,\gamma
\right) >0$\emph{\ depending only on listed parameters such that for any } $%
V^{\left( 1\right) },V^{\left( 2\right) }\in \overline{B\left( M\right)}$ 
\emph{the following estimate holds:} 
\begin{equation*}
\left\Vert J_{\lambda ,\gamma }^{\prime }\left( V^{\left( 2\right) }\right)
-J_{\lambda ,\gamma }^{\prime }\left( V^{\left( 1\right) }\right)
\right\Vert _{H_{2N}^{3}\left( \Omega \right) }\leq D\left\Vert V^{\left(
2\right) }-V^{\left( 1\right) }\right\Vert _{H_{2N}^{3}\left( \Omega
\right)}.
\end{equation*}

As to the existence and uniqueness of the minimizer, they are established in
Theorem 5.3. In fact, this theorem follows immediately from a combination of
above Theorems 5.1 and 5.2 with Lemma 2.1 and Theorem 2.1 of \cite{Bak}.
Hence, we omit its proof.

\textbf{Theorem 5.3}. \emph{Let the number }$\lambda _{2}=\lambda _{2}\left(
M,N,r,\Omega \right) >1$\emph{\ be the one in Theorem 5.1. Then for any }$%
\lambda \geq \lambda _{2}$\emph{\ and for any }$\gamma >0$\emph{\ the
functional }$J_{\lambda ,\gamma }\left( V\right) $\emph{\ has a unique
minimizer }$V_{\min ,\lambda ,\gamma }\in \overline{B\left( M\right) }$\emph{%
\ on }$\overline{B\left( M\right) }.$\emph{\ Furthermore, the following
inequality holds:}%
\begin{equation}
\left( J_{\lambda ,\gamma }^{\prime }\left( V_{\min ,\lambda ,\gamma
}\right) ,V_{\min ,\lambda ,\gamma }-Q\right) \leq 0\quad \text{for all }%
Q\in \overline{B\left( M\right) }.  \label{5.110}
\end{equation}

\subsection{The distance between the minimizer and the \textquotedblleft
ideal" solution}

\label{sec:5.3}

In accordance with the concept of Tikhonov for ill-posed problems \cite{T},
assume now that there exists the \textquotedblleft ideal" solution $V^{\ast
} $ of problem (\ref{eq:system})--(\ref{202}) with the \textquotedblleft
ideal" noiseless data $\psi _{0}^{\ast },\psi _{1}^{\ast }.$ It makes sense
to obtain an estimate of the distance between $V^{\ast }$ and the minimizer $%
V_{\min ,\lambda ,\gamma }$ of the functional $J_{\lambda ,\gamma }\left(
V\right) $ for the case of noisy data with the noise level $\delta \in
\left( 0,1\right) .$ This is what is done in the current subsection.

To obtain this estimate, we need to \textquotedblleft extend" the boundary
data $\psi _{0},\psi _{1}$ in (\ref{2}), (\ref{3}) inside $\Omega .$ Recall
that, unlike all previous works on the convexification, we have not done
this extension in the proof of our central Theorem 5.1. Thus, we assume
there exists a vector function $G\left( \mathbf{x}\right) \in
H_{2N}^{3}\left( \Omega \right) $ satisfying boundary conditions (\ref{2}), (%
\ref{3}),%
\begin{equation}
G\mid _{\partial \Omega }=\psi _{0}\left( \mathbf{x}\right) ,G_{z}\mathbf{%
\mid }_{\Gamma }=\psi _{1}\left( \mathbf{x}\right) .  \label{5.12}
\end{equation}%
On the other hand, the existence of the corresponding vector function $%
G^{\ast }\left( \mathbf{x}\right) \in H_{2N}^{3}\left( \Omega \right) $
satisfying boundary conditions with the \textquotedblleft ideal" data, 
\begin{equation}
G^{\ast }\mid _{\partial \Omega }=\psi _{0}^{\ast }\left( \mathbf{x}\right)
,G_{z}^{\ast }\mathbf{\mid }_{\Gamma }=\psi _{1}^{\ast }\left( \mathbf{x}%
\right)  \label{5.13}
\end{equation}%
follows from the existence of the ideal solution $V^{\ast }.$ We assume that 
\begin{equation}
\left\Vert G-G^{\ast }\right\Vert _{H_{2N}^{3}\left( \Omega \right) }<\delta
.  \label{5.14}
\end{equation}%
In addition, we suppose that%
\begin{equation}
\left\Vert V^{\ast }\right\Vert _{H_{2N}^{3}\left( \Omega \right)
},\left\Vert G^{\ast }\right\Vert _{H_{2N}^{3}\left( \Omega \right)
}<M-\delta .  \label{5.15}
\end{equation}%
Using (\ref{5.14}), (\ref{5.15}) and the triangle inequality, we easily see
that%
\begin{equation}
\left\Vert G\right\Vert _{H_{2N}^{3}\left( \Omega \right) }<M.  \label{5.16}
\end{equation}%
Our goal now is to estimate $\left\Vert V_{\min ,\lambda ,\gamma }-V^{\ast
}\right\Vert _{H_{2N}^{3}\left( \Omega \right) }$ via the noise parameter $%
\delta .$

\textbf{Theorem 5.4 }(accuracy and stability of minimizers).\emph{\ Suppose
that conditions (\ref{5.12})--(\ref{5.15}) hold. Let }$\lambda _{2}=\lambda
_{2}\left( M,N,r,\Omega \right) >1$\emph{\ be the number in Theorems 5.1,
5.3. Choose the number }$\lambda _{3}=\lambda _{2}\left( 3M,N,\Omega \right)
>\lambda _{2}>1$\emph{. Let }$\lambda =\lambda _{3}$\emph{\ and }$\gamma
=\delta ^{2}.$\emph{\ Then the following accuracy estimate holds}%
\begin{equation}
\left\Vert V_{\min ,\lambda ,\gamma }-V^{\ast }\right\Vert
_{H_{2N}^{2}\left( \Omega \right) }\leq C_{1}\delta .  \label{5.17}
\end{equation}

\textbf{Remark 5.1}. \emph{Since the power of }$\delta $\emph{\ is 1 in (\ref%
{5.17}), then it is natural to call (\ref{5.17}) \textquotedblleft the
Lipschitz stability estimate for the minimizers", which is similar to
Theorem 3.2. This estimate is obviously stronger than in all previous works
on the convexification; see, e.g. \cite%
{CAMWA,convexper,convIPnew,EIT,timedomain}, where one has }$\delta
^{\rho }$\emph{\ with }$\rho \in \left( 0,1\right) .$ \emph{The latter rate
is often called} \textquotedblleft \emph{the H\"{o}lder stability estimate
for the minimizers".}

\textbf{Proof of Theorem 5.4}. We note first that since the boundary
conditions for vector functions $V_{\min ,\lambda ,\gamma }$ and $V^{\ast }$
are different, then we cannot apply directly the strict convexity inequality
(\ref{5.1}) here, setting, e.g. that $V^{\left( 2\right) }=V^{\ast }$ and $%
V^{\left( 1\right) }=V_{\min ,\lambda ,\gamma }.$ And this is both the main
difficulty and the main new element of the proof, as compared with the
above-cited previous works on the convexification.

For every vector function $V\in B\left( M\right) $, consider the vector
function $W=V-G.$ Then by (\ref{5.16}) and the triangle inequality 
\begin{equation}
W\in B_{0}\left( 2M\right) =\left\{ W:\left\Vert W\right\Vert
_{H_{2N}^{3}\left( \Omega \right) }<2M,W\mid _{\partial \Omega }=W_{z}\mid
_{\Gamma }=0\right\} .  \label{5.180}
\end{equation}%
On the other hand, (\ref{5.16}) and (\ref{5.180}) imply that 
\begin{equation}
W+G\in B\left( 3M\right) \quad \text{for all }W\in B_{0}\left( 2M\right) .
\label{5.18}
\end{equation}%
Now, for any $W\in B_{0}\left( 2M\right) $ we have 
\begin{align}
& J_{\lambda ,\gamma }\left( W^{\ast }+G\right) -J_{\lambda ,\gamma }\left(
W+G\right) -J_{\lambda ,\gamma }^{\prime }\left( W+G\right) \left( W^{\ast
}-W\right)  \label{5.19} \\
& =J_{\lambda ,\gamma }\left( \widetilde{V}^{\ast }\right) -J_{\lambda
,\gamma }\left( V\right) -J_{\lambda ,\gamma }^{\prime }\left( V\right)
\left( \widetilde{V}^{\ast }-V\right) ,  \notag
\end{align}%
where $\widetilde{V}^{\ast }=W^{\ast }+G$ and $V=W+G.$ Notice that by (\ref%
{5.180}) and (\ref{5.18}) both vector functions $\widetilde{V}^{\ast },V\in
B\left( 3M\right) .$ Hence, by Theorem 5.1 we can apply the estimate (\ref%
{5.1}) to the second line of (\ref{5.19}) with $\lambda =\lambda
_{3}=\lambda _{2}\left( 3M,N,r,\Omega \right) >1.$ Thus, 
\begin{align}
& J_{\lambda _{3},\gamma }\left( W^{\ast }+G\right) -J_{\lambda _{3},\gamma
}\left( W+G\right) -J_{\lambda _{3},\gamma }^{\prime }\left( W+G\right)
\left( W^{\ast }-W\right)  \label{5.20} \\
& \geq C_{1}\left\Vert W^{\ast }-W\right\Vert _{H_{2N}^{2}\left( \Omega
\right) }^{2}+\gamma \left\Vert W^{\ast }-W\right\Vert _{H_{2N}^{3}\left(
\Omega \right) }^{2}\quad \text{for all }W\in B_{0}\left( 2M\right) .  \notag
\end{align}

Consider now the minimizer $V_{\min ,\lambda _{3},\gamma }\in \overline{%
B\left( M\right) }$ which is claimed by Theorem 5.3. Let $W_{\min ,\lambda
_{3},\gamma }=V_{\min ,\lambda _{3},\gamma }-G\in B\left( 2M\right) .$ Then (%
\ref{5.20}) implies that 
\begin{align}
& J_{\lambda _{3},\gamma }\left( W^{\ast }+G\right) -J_{\lambda _{3},\gamma
}\left( V_{\min ,\lambda _{3},\gamma }\right) -J_{\lambda _{3},\gamma
}^{\prime }\left( V_{\min ,\lambda _{3},\gamma }\right) \left( \left(
W^{\ast }+G\right) -V_{\min ,\lambda _{3},\gamma }\right)  \label{5.21} \\
& \geq C_{1}\left\Vert W^{\ast }-W_{\min ,\lambda _{3},\gamma }\right\Vert
_{H_{2N}^{2}\left( \Omega \right) }^{2}+\gamma \left\Vert W^{\ast }-W_{\min
,\lambda _{3},\gamma }\right\Vert _{H_{2N}^{3}\left( \Omega \right) }^{2}. 
\notag
\end{align}%
Using the triangle inequality, (\ref{5.14}) and (\ref{5.15}), we obtain%
\begin{align*}
& \left\Vert W^{\ast }+G\right\Vert _{H_{2N}^{3}\left( \Omega \right)
}=\left\Vert W^{\ast }+G^{\ast }+\left( G-G^{\ast }\right) \right\Vert
_{H_{2N}^{3}\left( \Omega \right) } \\
& \leq \left\Vert W^{\ast }+G^{\ast }\right\Vert _{H_{2N}^{3}\left( \Omega
\right) }+\left\Vert G-G^{\ast }\right\Vert _{H_{2N}^{3}\left( \Omega
\right) } \\
& =\left\Vert V^{\ast }\right\Vert _{H_{2N}^{3}\left( \Omega \right)
}+\left\Vert G-G^{\ast }\right\Vert _{H_{2N}^{3}\left( \Omega \right)
}<\left( M-\delta \right) +\delta =M.
\end{align*}%
This means that $\left( W^{\ast }+G\right) \in B\left( M\right) .$
Therefore, we use (\ref{5.110}) to get 
\begin{equation*}
-J_{\lambda _{3},\gamma }^{\prime }\left( V_{\min ,\lambda _{3},\gamma
}\right) \left( \left( W^{\ast }+G\right) -V_{\min ,\lambda _{3},\gamma
}\right) \leq 0.
\end{equation*}%
Hence, 
\begin{align*}
& J_{\lambda _{3},\gamma }\left( W^{\ast }+G\right) -J_{\lambda _{3},\gamma
}\left( V_{\min ,\lambda _{3},\gamma }\right) -J_{\lambda _{3},\gamma
}^{\prime }\left( V_{\min ,\lambda _{3},\gamma }\right) \left( \left(
W^{\ast }+G\right) -V_{\min ,\lambda _{3},\gamma }\right) \\
& \leq J_{\lambda _{3},\gamma }\left( W^{\ast }+G\right) .
\end{align*}%
Moreover, substituting this inequality in (\ref{5.21}), we obtain 
\begin{equation}
J_{\lambda _{3},\gamma }\left( W^{\ast }+G\right) \geq C_{1}\left\Vert
W^{\ast }-W_{\min ,\lambda _{3},\gamma }\right\Vert _{H_{2N}^{2}\left(
\Omega \right) }^{2}.  \label{5.22}
\end{equation}

We now estimate the left hand side of (\ref{5.22}). Note that the functional 
$J_{\lambda _{3},\gamma }\left( V\right) $ can be represented as 
\begin{align}
& J_{\lambda _{3},\gamma }\left( V\right) =J_{\lambda _{3},\gamma
}^{0}\left( V\right) +\gamma \left\Vert V\right\Vert _{H_{N}^{3}\left(
\Omega \right) }^{2},  \label{5.220} \\
& J_{\lambda _{3},\gamma }^{0}\left( V\right) =\exp \left[ -2\lambda \left(
R+r\right) ^{2}\right] \dint\limits_{\Omega }\left\vert L\left( V\right)
\right\vert ^{2}\mu _{\lambda }\left( z\right) d\mathbf{x}.  \label{5.221}
\end{align}%
Since $W^{\ast }+G^{\ast }=V^{\ast }$ is the ideal solution, then $L\left(
V^{\ast }\right) \left( \mathbf{x}\right) =0$ for $\mathbf{x}\in \Omega .$

Next, using the finite increment formula and (\ref{4.1}), we obtain%
\begin{align*}
\left\vert L\left( W^{\ast }+G\right) \right\vert ^{2}\left( \mathbf{x}%
\right) & =\left\vert L\left( W^{\ast }+G^{\ast }+G-G^{\ast }\right)
\right\vert ^{2}\left( \mathbf{x}\right) \\
& =\left\vert L\left( V^{\ast }\right) +S\left( G-G^{\ast }\right)
\right\vert ^{2}\left( \mathbf{x}\right) =\left\vert S\left( G-G^{\ast
}\right) \right\vert ^{2}\left( \mathbf{x}\right) ,
\end{align*}%
where by (\ref{40}) and (\ref{5.14}) the following estimate is valid:%
\begin{equation*}
\exp \left[ -2\lambda \left( R+r\right) ^{2}\right] \dint\limits_{\Omega
}\left\vert S\left( G-G^{\ast }\right) \right\vert ^{2}\left( \mathbf{x}%
\right) \mu _{\lambda }\left( z\right) d\mathbf{x}\leq C_{1}\delta ^{2}.
\end{equation*}%
This and (\ref{5.221}) imply that 
\begin{equation}
J_{\lambda _{3},\gamma }^{0}\left( W^{\ast }+G\right) \leq C_{1}\delta ^{2}.
\label{5.23}
\end{equation}%
Next, using (\ref{5.16}), (\ref{5.180}), (\ref{5.220}) and (\ref{5.23}), we
obtain 
\begin{equation}
J_{\lambda _{3},\gamma }\left( W^{\ast }+G\right) \leq C_{1}\left( \delta
^{2}+\gamma \right) .  \label{5.24}
\end{equation}%
Therefore, using (\ref{5.22}), (\ref{5.24}) and recalling that $\gamma
=\delta ^{2}$, we obtain%
\begin{equation}
\left\Vert W^{\ast }-W_{\min ,\lambda _{3},\gamma }\right\Vert
_{H_{2N}^{2}\left( \Omega \right) }\leq C_{1}\delta .  \label{5.26}
\end{equation}%
Finally, using (\ref{5.14}) and the triangle inequality, we obtain the
following lower bound for the left-hand side of (\ref{5.26})%
\begin{align*}
& \left\Vert W^{\ast }-W_{\min ,\lambda _{3},\gamma }\right\Vert
_{H_{2N}^{2}\left( \Omega \right) }=\left\Vert \left( W^{\ast }+G^{\ast
}\right) -\left( W_{\min ,\lambda _{3},\gamma }+G\right) +\left( G-G^{\ast
}\right) \right\Vert _{H_{2N}^{2}\left( \Omega \right) } \\
& =\left\Vert \left( V^{\ast }-V_{\min ,\lambda _{3},\gamma }\right) +\left(
G-G^{\ast }\right) \right\Vert _{H_{2N}^{2}\left( \Omega \right) } \\
& \geq \left\Vert V^{\ast }-V_{\min ,\lambda _{3},\gamma }\right\Vert
_{H_{2N}^{2}\left( \Omega \right) }-\left\Vert G-G^{\ast }\right\Vert
_{H_{2N}^{2}\left( \Omega \right) }\geq \left\Vert V^{\ast }-V_{\min
,\lambda _{3},\gamma }\right\Vert _{H_{2N}^{2}\left( \Omega \right) }-\delta
.
\end{align*}%
Substituting this in (\ref{5.26}), we obtain the target estimate (\ref{5.17}%
). $\ \ \ \ \square $

\textbf{Corollary 5.1}. \emph{The functional }$I_{\lambda ,\gamma }\left(
W\right) :=J_{\lambda ,\gamma }$\emph{\ }$\left( W+G\right) $\emph{\ is
strictly convex on }$\overline{B_{0}\left( 2M\right) }$\emph{\ for all }$%
\lambda \geq \lambda _{3},$\emph{\ where }$\lambda _{3}$\emph{\ is the
number defined in Theorem 5.4.}

\textbf{Proof}. It follows from the proof of Theorem 5.4 and (\ref{5.1})
that the following analog of (\ref{5.20}) holds for all $\lambda \geq
\lambda _{3}$ and for all $W^{\left( 1\right) },W^{\left( 2\right) }\in
B_{0}\left( 2M\right)$%
\begin{align*}
& I_{\lambda ,\gamma }\left( W^{\left( 2\right) }\right) -I_{\lambda ,\gamma
}\left( W^{\left( 1\right) }\right) -I_{\lambda ,\gamma }^{\prime }\left(
W^{\left( 1\right) }\right) \left( W^{\left( 2\right) }-W^{\left( 1\right)
}\right) \\
& \geq C_{1}\left\Vert W^{\left( 2\right) }-W^{\left( 1\right) }\right\Vert
_{H_{2N}^{2}\left( \Omega \right) }^{2}+\gamma \left\Vert W^{\left( 2\right)
}-W^{\left( 1\right) }\right\Vert _{H_{2N}^{3}\left( \Omega \right) }^{2}. \
\ \ \ \ \ \ \ \ \ \ \ \square
\end{align*}

\section{The Globally Convergent Gradient Projection Method}

\label{sec:6}

Now we construct an approximation for the vector function $W^{\ast }=V^{\ast
}-G^{\ast }$ for $W^{\ast }\in B_{0}\left( 2M\right) .$ It follows from (\ref%
{5.180}) that $B_{0}\left( 2M\right) \subset H_{0,2N}^{3}\left( \Omega
\right) .$ Let $P_{\overline{B}}:H_{0,2N}^{3}\left( \Omega \right)
\rightarrow \overline{B_{0}\left( 2M\right) }$ be the orthogonal projection
operator of the space $H_{0,2N}^{3}\left( \Omega \right) $ on the closed
ball $\overline{B_{0}\left( 2M\right) }.$ Let $W^{\left( 0\right) }\in
B_{0}\left( 2M\right) $ be an arbitrary point of the ball $B_{0}\left(
2M\right) $ and let $\eta >0$ be a number. The gradient projection method
constructs the following sequence:%
\begin{equation}
W^{\left( n\right) }=P_{\overline{B}}\left( W^{\left( n-1\right) }-\eta
J_{\lambda ,\gamma }^{\prime }\left( W^{\left( n-1\right) }+G\right) \right)
,\text{ }n=1,2,...  \label{6.1}
\end{equation}%
It is important for computations that $\left( W^{\left( n-1\right) }-\eta
J_{\lambda ,\gamma }^{\prime }\left( W^{\left( n-1\right) }+G\right) \right)
=:Y_{n-1}\left( \mathbf{x}\right) \in H_{0,2N}^{3}\left( \Omega \right) .$
Indeed, $W^{\left( n-1\right) }\in H_{0,2N}^{3}\left( \Omega \right) $;
also, (\ref{555}) holds. In other words, the vector function $%
Y_{n-1}\left( \mathbf{x}\right) $ satisfies the boundary conditions $%
Y_{n-1}\mid _{\partial \Omega }=\partial _{z}Y_{n-1}\mid _{\Gamma }=0.$

\textbf{Theorem 6.1}. \emph{Let }$\lambda \geq \lambda _{2}^{\prime
}=\lambda _{2}\left( 2M,N,r,\Omega \right) >1,$\emph{\ where }$\lambda _{2}$%
\emph{\ was defined in Theorem 5.1. Let }$W_{\min ,\lambda ,\gamma }$\emph{\
be the minimizer of the functional }$J_{\lambda ,\gamma }\left( W+G\right) $%
\emph{\ on the set }$\overline{B_{0}\left( 2M\right) }$\emph{, the existence
and uniqueness of which follow from Theorem 5.3 and Corollary 5.1. Then
there exists a sufficiently small number }$\eta _{0}=\eta _{0}\left(
2M,N,r,\Omega ,\lambda \right) \in \left( 0,1\right) $\emph{\ depending only
on listed parameters such that for any }$\eta \in \left( 0,\eta _{0}\right) $%
\emph{\ we can find a number }$\theta =\theta \left( \eta \right) \in \left(
0,1\right) $ \emph{such that the sequence }$\left\{ W^{\left( n\right)
}\right\} _{n\in \mathbb{N}^{\ast }}$\emph{\ converges to }$W_{\min ,\lambda
,\gamma }$\emph{\ in the }$H_{2N}^{3}\left( \Omega \right) $--\emph{norm}%
\emph{\ and the following convergence estimate holds:}%
\begin{equation}
\left\Vert W_{\min ,\lambda ,\gamma }-W^{\left( n\right) }\right\Vert
_{H_{2N}^{3}\left( \Omega \right) }\leq \theta ^{n}\left\Vert W_{\min
,\lambda ,\gamma }-W^{\left( 0\right) }\right\Vert _{H_{2N}^{3}\left( \Omega
\right) },\;n=1,2,...  \label{6.2}
\end{equation}

Theorem 6.1 follows immediately from the combination of Theorems 5.1--5.3
with Theorem 2.1 of \cite{Bak}.

\textbf{Theorem 6.2}. \emph{Let }$\lambda =\lambda _{4}=\lambda _{2}\left(
2M,N,r,\Omega \right) \geq \lambda _{2}^{\prime }.$\emph{\ Suppose that
conditions imposed in Theorems 5.4 and 6.1 hold. Then the following
convergence estimates are valid for }$n=1,2,...$\emph{\ }%
\begin{align}
& \left\Vert W^{\ast }-W^{\left( n\right) }\right\Vert _{H_{2N}^{2}\left(
\Omega \right) }\leq C_{1}\delta +\theta ^{n}\left\Vert W_{\min ,\lambda
,\gamma }-W^{\left( 0\right) }\right\Vert _{H_{2N}^{3}\left( \Omega \right)
},  \label{6.3} \\
& \left\Vert c^{\ast }\left( \mathbf{x}\right) -c_{n}\left( \mathbf{x}%
\right) \right\Vert _{L^{2}\left( \Omega \right) }\leq C_{1}\delta +\theta
^{n}\left\Vert W_{\min ,\lambda ,\gamma }-W^{\left( 0\right) }\right\Vert
_{H_{2N}^{3}\left( \Omega \right) },  \label{6.4}
\end{align}%
\emph{where }$c^{\ast }\left( \mathbf{x}\right) $\emph{\ stands in the
right hand side of equation (\ref{3.02}) in the case when }$W^{\ast }\left( 
\mathbf{x}\right) $\emph{\ is replaced with }$V^{\ast }\left( \mathbf{x}%
\right) =W^{\ast }\left( \mathbf{x}\right) +G^{\ast }\left( \mathbf{x}%
\right) $\emph{. The function }$v^{\ast }\left( \mathbf{x},\alpha \right) $%
\emph{\ is obtained via components of the vector function }$V^{\ast }\left( 
\mathbf{x}\right) $\emph{\ and (\ref{eq:truncated}) and then this function
is substituted in the left hand side of (\ref{3.03}); see the first item of
Remarks 3.1. The function }$c_{n}\left( \mathbf{x}\right) $\emph{\ is
obtained in the same way with the only replacement of }$V^{\ast }\left( 
\mathbf{x}\right) $\emph{\ with }$V^{\left( n\right) }\left( \mathbf{x}%
\right) =W^{\left( n\right) }\left( \mathbf{x}\right) +G\left( \mathbf{x}%
\right) .$

\textbf{Proof}. Combining (\ref{5.26}) with (\ref{6.2}) we obtain%
\begin{align*}
& \left\Vert W^{\ast }-W^{\left( n\right) }\right\Vert _{H_{2N}^{2}\left(
\Omega \right) }=\left\Vert \left( W^{\ast }-W_{\min ,\lambda ,\gamma
}\right) +\left( W_{\min ,\lambda ,\gamma }-W^{\left( n\right) }\right)
\right\Vert _{H_{2N}^{2}\left( \Omega \right) } \\
& \leq \left\Vert W^{\ast }-W_{\min ,\lambda ,\gamma }\right\Vert
_{H_{2N}^{2}\left( \Omega \right) }+\left\Vert W_{\min ,\lambda ,\gamma
}-W^{\left( n\right) }\right\Vert _{H_{2N}^{2}\left( \Omega \right) } \\
& \leq C_{1}\delta +\theta ^{n}\left\Vert W_{\min ,\lambda ,\gamma
}-W^{\left( 0\right) }\right\Vert _{H_{2N}^{3}\left( \Omega \right) },
\end{align*}%
which proves (\ref{6.3}). As to (\ref{6.4}), it follows from (\ref{6.3}) and
the construction of functions $c^{\ast }\left( \mathbf{x}\right)
,c_{n}\left( \mathbf{x}\right) $ described in the formulation of Theorem
6.2. $\ \ \ \ \square $

\textbf{Remarks 6.1}.

\begin{enumerate}
\item \emph{Since the starting point }$W^{\left( 0\right) }$\emph{\ of the
gradient projection method (\ref{6.1}) is an arbitrary point of the ball }$%
B_{0}\left( 2M\right) $\emph{\ and since smallness conditions are not
imposed on }$M$\emph{, then convergence estimates (\ref{6.3}) and (\ref{6.4}%
) mean the \textbf{global convergence} of the gradient projection method (%
\ref{6.1}) to the correct solution. In other words, a good first guess about
the ideal solution is no longer required. We note that in the case of a non
convex functional, the global convergence of a gradient-like method cannot
be guaranteed.}

\item \emph{Another observation here is that due to estimate (\ref{5.17}) we
have \textquotedblleft Lipschitz-like" convergence rate in (\ref{6.3}), (\ref%
{6.4}) with respect to }$\delta $\emph{. This is stronger than the
\textquotedblleft H\"{o}lder-like" convergence rates in all previous
publications \cite{CAMWA,convexper,convIPnew,EIT,timedomain} about
the convexification where }$\delta $\emph{\ is replaced with }$\delta ^{\rho
},$\emph{\ where }$\rho \in \left( 0,1\right) ,$ \emph{also, see Remark 5.1.}
\end{enumerate}

\section{Numerical Results}

\label{sec:7}

\subsection{Some details of the numerical implementation}

\label{sec:7.1}

First, we note that it follows from (\ref{300}), (\ref{eq:truncated}) and (%
\ref{301}) that%
\begin{equation}
V\mid _{\partial \Omega \diagdown \Gamma }=0.  \label{7.1}
\end{equation}%
Let $\left\{ \alpha _{n}\right\} _{n=0}^{\ell -1}$ be the set of selected
ode points on the segment $\left[ -a,a\right] ,$ i.e. $-a=\alpha _{0}<\alpha
_{1}<\cdots <\alpha _{\ell -1}=a,\text{ }\alpha _{n}-\alpha _{n-1}=h_{s},$
where $h_{s}>0$ is the grid step size. Thus, in our computations, our
sources are $\left\{ \mathbf{x}_{\alpha _{n}}\right\} _{n=0}^{l-1}=\left\{
\left( \alpha _{n},0,-d\right) \right\} _{n=0}^{l-1}\subset L_{\text{src}};$
see (\ref{2001}). To minimize the functional $J_{\lambda ,\gamma }\left(
V\right) $ in (\ref{4.2}), we rewrite the involved derivatives via finite
differences. The idea is to minimize the resulting functional with respect
to the values of $V$ at grid points. We use the same grid step size $h$ in $%
x,y$,$z$ directions. For any vector function $u\left( \mathbf{x}\right) $ we
denote $u_{p,q,s}=u(x_{p},y_{q},z_{s})$ the corresponding discrete function
defined at grid points.

We define the Laplace operator in finite differences as $\Delta
^{h}u_{p,q,s}=\partial _{xx}^{h}u_{p,q,s}+\partial
_{yy}^{h}u_{p,q,s}+\partial _{zz}^{h}u_{p,q,s}$, where, for interior grid
points of $\Omega $ we use, e.g., 
\begin{align*}
\partial _{xx}^{h}u_{p,q,s}=h^{-2}\left(
u_{p+1,q,s}-2u_{p,q,s}+u_{p-1,q,s}\right) .
\end{align*}
Then, the gradient operator in finite difference $\nabla
^{h}u_{p,q,s}=\left( \partial _{x}^{h}u_{p,q,s},\partial
_{y}^{h}u_{p,q,s},\partial _{z}^{h}u_{p,q,s}\right) $ follows. In addition,
the data at $\Gamma$ are given by $\partial _{z}^{h}u_{p,q,0}=h^{-1}\left(
u_{p,q,1}-u_{p,q,0}\right) . $

We have applied the matrix $S_{N}^{-1}$ in (\ref{202}) to obtain equation (%
\ref{eq:system}). However, this is convenient only for the above theory. In
computations we do not apply $S_{N}^{-1}.$ The resulting matrix equation is
equivalent to (\ref{eq:system}) and analogs of the above Theorems 5.2-5.4,
6.1, 6.2 can be straightforwardly formulated for the continuous form of
functional (\ref{eq:discreteJ}), which is the direct analog of functional (%
\ref{4.2}). Now about the computational implementation of boundary
conditions. Rather than satisfying boundary conditions exactly, we minimize
the differences between boundary values of the discrete vector function $%
V_{p,q,s}$ and boundary conditions. To do this, we use penalty terms with
certain weights $K_{0},K_{1},K_{2}>0$. These weights are chosen numerically.

Let $V^{h}=\left\{ V_{p,q,s}^{h}\right\} _{p,q,s=0}^{Z_{h}-1}$ be the
discrete version of the vector function $V$. Thus, taking into account (\ref%
{eq:system})--(\ref{3}) and (\ref{7.1}), for each $N\geq 1$ we minimize the
following fully discrete form of the weighted Tikhonov-like $J_{\lambda
,\gamma }\left( V\right) $:

\begin{align}
&J_{\lambda ,\gamma }^{h}\left( V^{h}\right) =e^{-2\lambda \left( R+r\right)
^{2}} \sum_{p,q,s=0}^{Z_{h}-1}h^{3}\left\vert S_{N}\Delta
^{h}V_{p,q,s}^{h}+f\left( \nabla ^{h}V_{p,q,s}^{h}\right) \right\vert
^{2}\mu _{\lambda }\left( z_{s}\right)  \label{eq:discreteJ} \\
& +\sum_{p,q=0}^{Z_{h}-1}h^{2}\left( K_{0}\left\vert V_{p,q,0}^{h}-\psi
_{0,p,q}\right\vert ^{2}+K_{1}\left\vert \partial _{z}^{h}V_{p,q,0}^{h}-\psi
_{1,p,q}\right\vert ^{2}\right)  \notag \\
& +K_{2}\sum_{p,q=0}^{Z_{h}-1}h^{2}\left\vert w_{p,q,Z_{h}-1}\right\vert
^{2}+K_{2}\sum_{q,s=0}^{Z_{h}-1}h^{2}\left( \left\vert w_{0,q,s}\right\vert
^{2}+\left\vert w_{Z_{h}-1,q,s}\right\vert ^{2}\right)  \notag \\
& +K_{2}\sum_{p,s=0}^{Z_{h}-1}h^{2}\left( \left\vert w_{p,0,s}\right\vert
^{2}+\left\vert w_{p,Z_{h}-1,s}\right\vert ^{2}\right) +\gamma
\sum_{p,q,s=0}^{Z_{h}-1}h^{3}\left( \left\vert V_{p,q,s}^{h}\right\vert
^{2}+\left\vert \nabla ^{h}V_{p,q,s}^{h}\right\vert ^{2}\right) .  \notag
\end{align}

It follows from (\ref{eq:discreteJ}) that we consider the regularization
term in the $H^{1}$ norm. Indeed, it is much more complicated to implement
the $H^{3}$ norm. On the other hand, if the number of grid points is not too
large, then these norms, taken in the discrete forms, are \textquotedblleft
effectively" equivalent. Our numerical experience shows that the $H^{1}$
norm in the regularization term is sufficient. Terms with $%
K_{0},K_{1},K_{2}>0$ are those penalty terms mentioned above with respect to
the boundary conditions. Since these terms are convex, then they do not ruin
the strict convexity of our functional. Numbers $K_{0},K_{1},K_{2}$ are
chosen numerically and will be specified in \cref{sec:experiments}.


\subsection{Numerical studies}

\label{sec:main}

\subsubsection{Generic algorithm}

\label{sec:alg}

For each $\alpha \in \lbrack -a,a]$, we computationally simulate the data (%
\ref{303}) by solving the Lippmann--Schwinger equation (\ref{302}). Thus, on
the grid of $\mathbf{x}\in \Gamma $ and $\mathbf{x}_{\alpha }\in L_{\text{src%
}}$ one has $u_{p,q,s}=F(x_{p},y_{q},z_{s},\mathbf{x}_{\alpha _{l}})$ for
each $l=\overline{0,\ell -1}$. Given $N$ in (\ref{eq:truncated}) and $k>0$,
our analysis then leads to the following algorithm:

\begin{enumerate}
\item To generate the data $\psi _{j}(x_{p},y_{q},z_{s})$ for $j=0,1$ in (%
\ref{2}), (\ref{3}), solve the Lippmann--Schwinger equation (\ref{302}).

\item Compute the minimizer $V_{\min }^{h}$ of the functional $J_{\lambda
,\gamma }(V^{h})$ in \eqref{eq:discreteJ} as an approximation of $V^{h}$ .

\item Using the minimizer of item 2 and the special basis $\left\{ \Psi
_{n}(\alpha )\right\} _{n=0}^{N-1}$, construct an approximation of the
auxiliary function $v(x_{p},y_{q},z_{s},\alpha _{l})$;

\item Compute an approximation of the unknown dielectric constant $%
c(x_{p},y_{q},z_{s},k)$ by the following formulae: 
\begin{equation}  \label{ccc}
c_{p,q,s}=\underset{\alpha _{l}}{\text{mean}}\left\vert -\frac{\Delta
^{h}v_{p,q,s,\alpha _{l}}+\left( \nabla ^{h}v_{p,q,s,\alpha _{l}}\right)
^{2}+2\nabla ^{h}v_{p,q,s,\alpha _{l}}\cdot \tilde{\mathbf{x}}_{p,q,s,\alpha
_{l}}}{k^{2}}\right\vert +1,
\end{equation}%
where $\tilde{\mathbf{x}}_{p,q,s,\alpha _{l}}$ denotes the value of $\tilde{%
\mathbf{x}}_{\alpha }$ at $(x_{p},y_{q},z_{s})$ for every $\alpha _{l}$; see
(\ref{400}).
\end{enumerate}

We use the absolute value here since one should have $c_{p,q,s}\geq 1,$ see (%
\ref{2020}).

\begin{figure}[]
\centering
\begin{subfigure}{0.42\linewidth}
		\includegraphics[width=\linewidth]{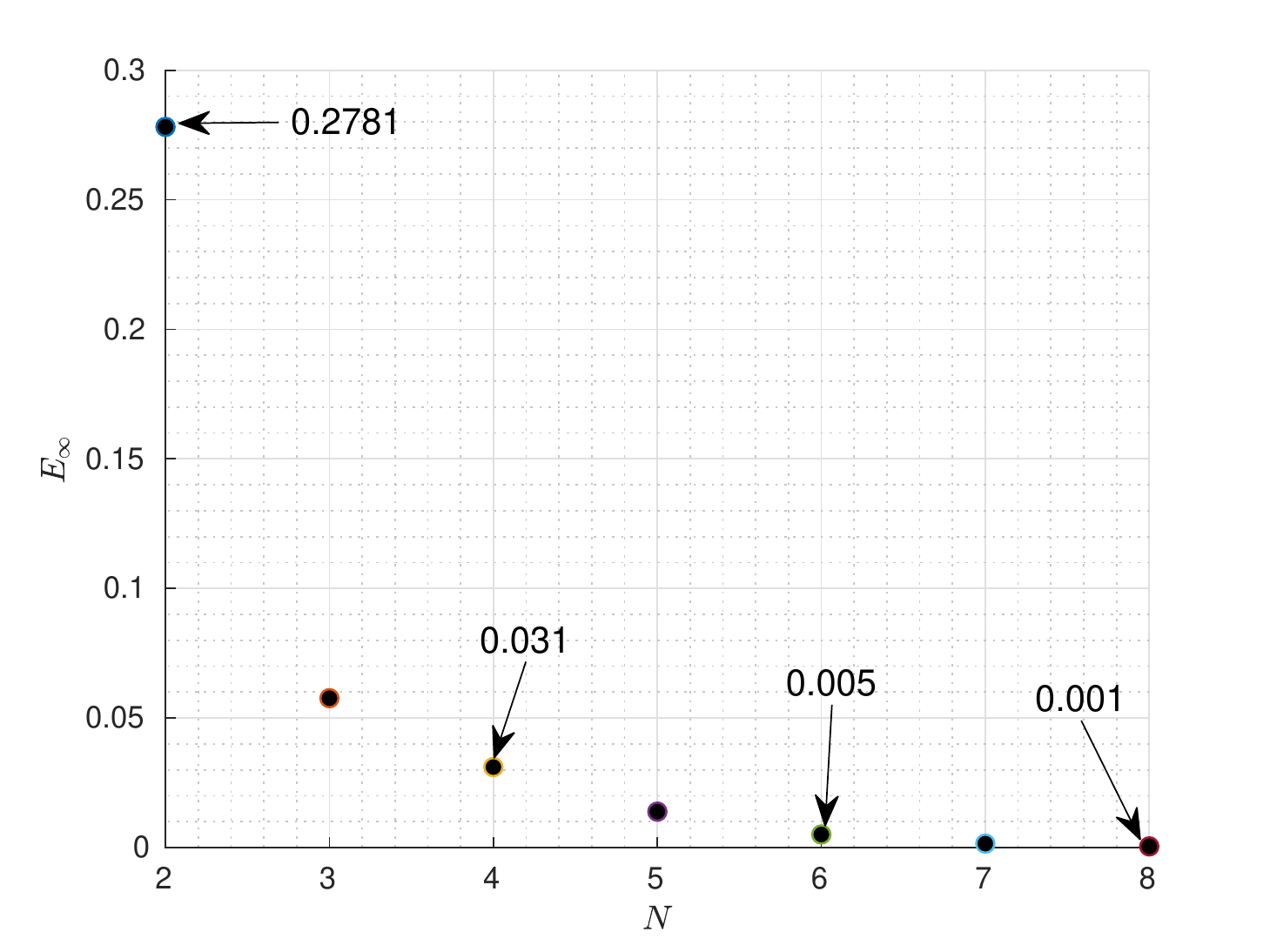}
		\caption{Relative $L^{\infty}$ error\label{fig:0a}}
	\end{subfigure}
\begin{subfigure}{0.42\linewidth}
		\includegraphics[width=\linewidth]{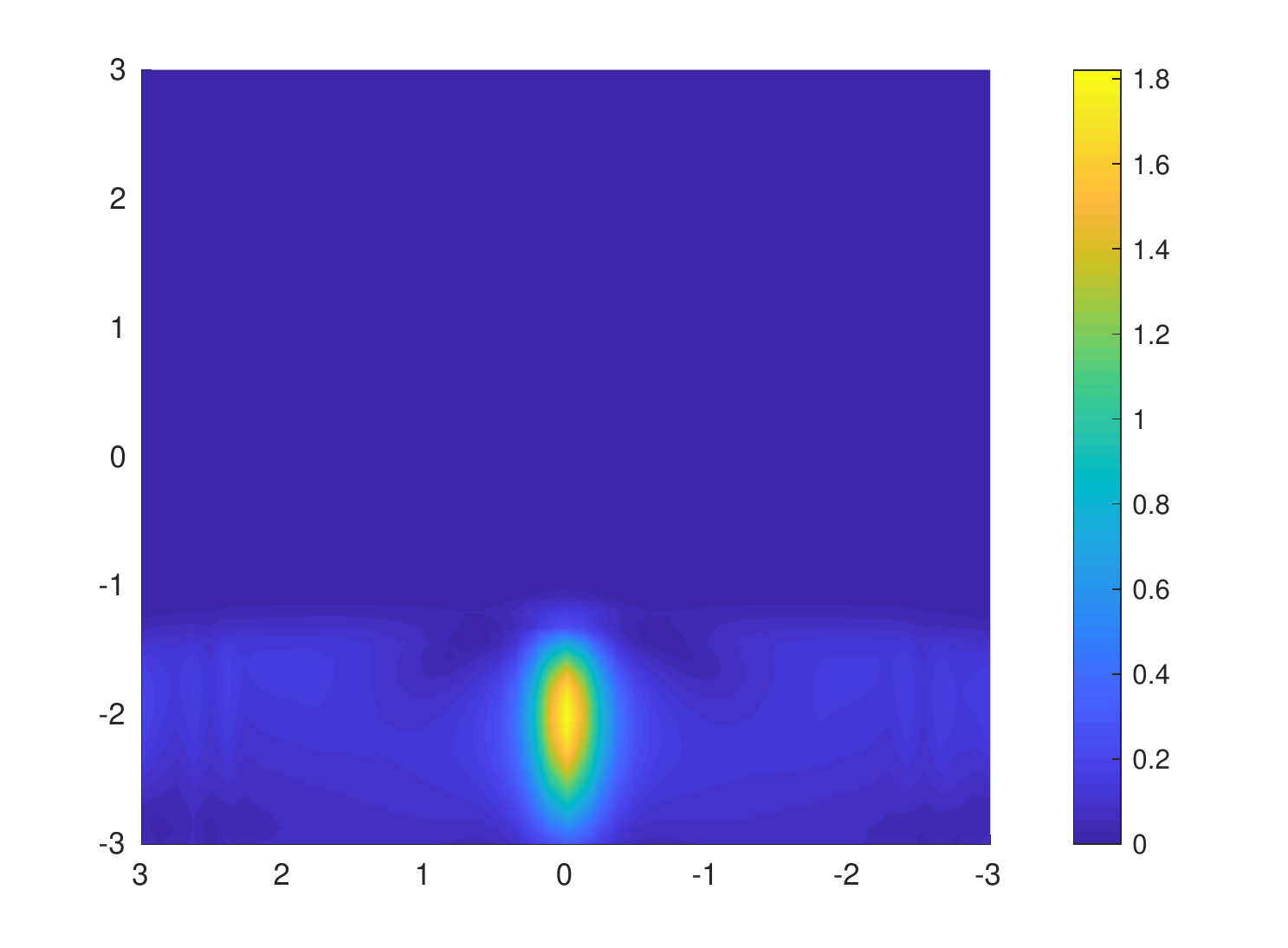}
		\caption{Cross-section after the first step\label{fig:0f}}
	\end{subfigure}
\caption{(a) The choice of $N$ based on the relative $L^{\infty}$ error. (b)
Cross-section by the plane $\left\{x=0\right\}$ of the image computed in
Test 2 after finishing the first step of computations.}
\label{fig:0}
\end{figure}

\subsubsection{Computational procedure}

\label{sec:experiments}

We now describe some details of the algorithm of section \cref{sec:alg}. The
domain under consideration is now the three-dimensional cube with the edge
length $R=3$. Besides, we choose $a=1$, $\gamma =10^{-4}$, $d=7.5$, $N=4$, $%
Z_{h}=51$ and $\ell =11$. Also, we choose $K_{0}=1,K_{1}=2$ and $%
K_{2}=10^{-3}$ in \eqref{eq:discreteJ} for all steps. The wavenumber $k$ is
specified below. Even though our above theory says that we need to use the
gradient projection method and large values of $\lambda $, we have
discovered computationally that the simpler to implement gradient descent
method and a reasonable value $\lambda =1.1$ work well. So, we use them.
These observations coincide with those of our previous works on the
numerical studies of the convexification \cite%
{CAMWA,convexper,convIPnew,EIT,timedomain}.

The choice of $N$ can be specified as follows. We take a reference sample to
be imaged. This is the one of Test 1; see \cref{table:1} for details. Once
chosen, we keep $N$ the same for all other examples. First, we solve the
Lippmann--Schwinger equation \eqref{302} to generate the data for the
reference inclusion of Test 1. Then we obtain $u_{\text{true}}\left( \mathbf{%
x},\alpha \right) $ and $v_{\text{true}}\left( \mathbf{x},\alpha \right) $,
respectively. Denote $v_{\text{true},n}\left( \mathbf{x}\right) $ the
corresponding Fourier coefficient of $v_{\text{true}}\left( \mathbf{x}%
,\alpha \right) $ with respect to our basis $\left\{ \Psi _{n}(\alpha
)\right\} _{n=0}^{\infty }$. Having $v_{\text{true},n}(\mathbf{x})$
numerically, we can compute the function $v_{\text{true}}^{N}(\mathbf{x,}%
\alpha )$ in (\ref{eq:truncated}). Hence, we are able to compute the
following relative $L^{\infty }$-like error: 
\begin{equation*}
E_{\infty }\left( v_{\text{true}}\right) =\frac{{\displaystyle\max_{\mathbf{x%
},\alpha }\left\vert v_{\text{true}}-v_{\text{true}}^{N}\right\vert }}{{%
\displaystyle\max_{\mathbf{x},\alpha }\left\vert v_{\text{true}}\right\vert }%
}\times 100\%.
\end{equation*}%
We observe in \cref{fig:0a} that $N=4$ is acceptable in the sense that the
error $E_{\infty }\left( v_{\text{true}}\right) $ is sufficiently small
(around $5\times 10^{-2}$). Hence, it is intuitively clear that increasing $%
N $ would only increase the computational time without providing an
essential difference in results. This choice of $N$ is in an agreement with,
e.g., \cite{Khoa} where a truncation was done to solve an inverse boundary
value problems of an elliptic equation.

As the number of point sources may not be large in practice, we compute $v_{%
\text{true},n}$ by using the Gauss--Legendre quadrature method: 
\begin{equation*}
v_{\text{true},n}(\mathbf{x}) := \left\langle v_{\text{true}}\left(\mathbf{x}%
,\cdot\right),\Psi_{n}\left(\cdot\right)\right\rangle_{L^2(-a,a)} \approx
\sum_{l=0}^{\ell}\text{w}_{l}v_{\text{true}}\left(\mathbf{x},\tilde{\alpha}%
_{l}\right)\Psi_{n}(\tilde{\alpha}_{l}),
\end{equation*}
where $\tilde{\alpha}_{l}$ are the abscissae in $[-a,a]$ and $\text{w}_{l}$
are the corresponding weights. Since these abscissae are fixed, we get the
values of $v_{\text{true}}\left(\mathbf{x},\tilde{\alpha}_{l}\right)$ from $%
v_{\text{true}}\left(\mathbf{x},\alpha_{l}\right)$ using the spline
interpolation via the built-in \texttt{griddedInterpolant} in \textsc{Matlab}%
. As we compute the basis $\Psi_{n}$ symbolically, we know their values at
the abscissae very precisely.


Following the algorithm of \cref{sec:alg}, we now detail the computational
procedure in the minimization process applied to the functional (\ref%
{eq:discreteJ}). Recall that the Dirichlet data at $\partial \Omega
\diagdown \Gamma $ are obtained by the heuristic data completion (\ref{301})
and this completion forms the final item of our approximate mathematical
model (the first item of Remarks 3.1). Since the magnitude of the
backscattering data is small and since the correct boundary data on $%
\partial \Omega \diagdown \Gamma $ are actually neglected by (\ref{301}),
then one can anticipate that only those targets will be reasonably imaged,
which are close to the measurement side $\Gamma $ of the cube $\Omega $. In
other words, the data propagation procedure mentioned in \cref{sec:3.3}
should estimate well distances to targets, and it was shown in section 6.2
of \cite{Nguyen2018} for the case of experimentally collected data that the
latter is possible. It is also clear that this very limited information
indicates that our reconstruction is very challenging.

On the other hand, it follows from (\ref{4.03}) and Theorem 5.1 that all
iterates of the gradient descent method we use should have the same
Dirichlet boundary conditions at $\partial \Omega $ and the same Neumann
boundary condition at $\Gamma $ (i.e. at $z=-R$)$.$ Furthermore, our
computational experience shows that since the correct boundary data on $%
\partial \Omega \diagdown \Gamma $ are not given, then we need to have a
partly heuristic method of estimating the location and sizes of the target
to be imaged. These cause our choice of the starting point of iterations as
well as the choice of first and second steps of our numerical procedure
described below. Nevertheless, the main point is that our choices still use
only the measured data and do not rely on any information about a small
neighborhood of the correct solution.

\emph{Step 1.} This step consists of two (2) substeps. As it was pointed out
in the Introduction section, the focus of our application is the detection
and identification of antipersonnel land mines and IEDs. The sizes of these
targets are usually small: between 5 and 15 centimeters (cm). Therefore, we
search for targets in the domain $\Omega_{1}:=\left\{-R \le z \le -R +
2\right\}\subset \Omega$, which means 20 cm in depth from the ground
boundary. In the first substep, we choose the following linear
approximation, denoted by $V_{0}\left( \mathbf{x}\right) \in \mathbb{R}^{N}$%
, as the starting point of iterations in the minimization of the functional $%
J_{\lambda ,\gamma }^{h}\left( V^{h}\right) $: 
\begin{align}  \label{304}
& V_{0}\left( \mathbf{x}\right) =%
\begin{pmatrix}
v_{00}\left( \mathbf{x}\right), & v_{01}\left( \mathbf{x}\right), & \cdots,
& v_{0(N-1)}\left( \mathbf{x}\right)%
\end{pmatrix}%
^{T}, \\
& v_{0n}=\left( \psi _{0n}+\psi _{1n}(z+R)\right) \chi (z),  \label{3042}
\end{align}%
where $\chi :[-R,R]\rightarrow \mathbb{R}$ is a smooth function given by 
\begin{equation}
\chi\left(z\right)=%
\begin{cases}
\exp\left(\frac{\left(z+R\right)^{2}}{\left(z+R\right)^{2}-\left(R-1%
\right)^{2}}\right) & \text{if }z<-1, \\ 
0 & \text{otherwise.}%
\end{cases}
\label{1000}
\end{equation}
It is worth noting that $\chi$ attains the maximum value of 1 at $z=-R$ and
then, we can show that $v_{0n}\left(z=-R\right) = \psi_{0n}$, $%
\partial_{z}v_{0n}\left(z=-R\right) = \psi_{1n}$. Thus, the starting point $%
v_{0n}$ satisfies the Dirichlet and Neumann boundary conditions at $z=-R$.
Besides, the vector function $V_0$ satisfies the zero Dirichlet boundary
condition at $z=R$ because $\chi$ tends to 0 as $z\to -1^{+}$. We have
numerically observed that $\psi _{0n}=\psi _{1n}=0$ at $x,y=\pm R.$ Hence,
the vector function defined in (\ref{304}) satisfies the zero Dirichlet
boundary conditions. On the first substep, we minimize the functional $%
J_{\lambda ,\gamma }^{h}\left( V^{h}\right) $ with the starting point (\ref%
{304}).

As to the step size $\eta $ of the gradient descent method, we start from $%
\eta _{1}=10^{-1}.$ This $\eta _{1}$ is unchanged on the first step as well
as the second step below. On each iterative step number $m$ the step size $%
\eta _{m}$ is reduced by the factor of 2 if the value of the functional at
the iteration number $m$ exceeds its value on the previous iteration, i.e.
if $J_{\lambda ,\gamma }^{h}\left( V_{m}^{h}\right) \geq J_{\lambda ,\gamma
}^{h}\left( V_{m-1}^{h}\right) $. Otherwise, $\eta _{m+1}=\eta _{m}.$ The
minimization process is stopped at $m_{\text{stop}}$ when either $\eta _{m_{%
\text{stop}}}<10^{-8}$ or $\left\vert J_{\lambda ,\gamma }^{h}\left( V_{m_{%
\text{stop}}}^{h}\right) -J_{\lambda ,\gamma }^{h}\left( V_{m_{\text{stop}%
}-1}^{h}\right) \right\vert <10^{-8}.$




After the first substep, we obtain the numerical coefficient of $c\left( 
\mathbf{\ x}\right) $, denoted by $\tilde{c}\left( \mathbf{\ x}\right)$.
Then we complete the first step of computations by getting rid of possible
artifacts in $\tilde{c}$ via a postprocessing, which forms the second
substep. This is done by replacing the function $\tilde{c}$ with the
function $c_{\text{temp}},$ where 
\begin{equation*}
c_{\text{temp}}\left(\mathbf{x}\right)=%
\begin{cases}
\tilde{c}\left(\mathbf{x}\right) & \text{if }\left|\hat{c}\left(\mathbf{x}%
\right)\right|\ge0.7\max_{\mathbf{x}}\left|\tilde{c}\left(\mathbf{x}%
\right)\right|, \\ 
0 & \text{otherwise.}%
\end{cases}%
\end{equation*}
Lastly, the function $c_{\text{temp}}$ is smoothed by a Gaussian filter. We
will shortly specify the Gaussian procedure in the end of Step 2. 

\emph{Step 2.} We have numerically observed that, in the first step, we can
find a good approximation for the location of the target and a somewhat good
approximation of its shape. The second step is for an improvement of the
values of the function $c\left( \mathbf{x}\right) $. We start this step by
plugging $c_{\text{temp}}$ in the Lippmann--Schwinger equation (\ref{302}),
solving it, and thus, obtaining a new vector function $\hat{v}\left( \mathbf{%
x},\alpha \right) $. Then we find first $N$ Fourier components of $\hat{v}%
\left( \mathbf{x},\alpha \right) $, denoted by $\hat{v}_{n}\left( \mathbf{x}%
\right) $, as in (\ref{eq:truncated}). Note that after Step 1, we can see
where the object is located by an evaluation of 2D cross-sections of the
image obtained after Step 1; see e.g. \cref{fig:0f}. This way enables us to
narrow the domain of our search, i.e. we look for the object in $\Omega
_{2}:=\left\{ -b_{x}\leq x\leq b_{x},-b_{y}\leq y\leq b_{y},-R\leq z\leq
-b_{z}\right\} \subset \Omega _{1},$ for some numbers $b_{x},b_{y},b_{z}>0$.
We rely on this information to get a smooth function, denoted by $\tilde{\chi%
}(\mathbf{x})$, with $\hat{v}_{n}$ and ensure the boundary conditions during
the minimization. Recall the function $\chi $ in (\ref{1000}) which is
essentially generated by the function $\exp \left( \frac{t^{2}}{t^{2}-1}%
\right) $ for $t\in \lbrack 0,1)$ tending to 0 as $t\rightarrow 1^{+}$ and
attaining the maximum 1 at $t=0$. Moreover, the derivative of this function
vanishes at $t=0$. Thereby, the starting point of iterations of the gradient
descent method for the second stage of the minimization of the functional $%
J_{\lambda ,\gamma }^{h}\left( V^{h}\right) $ is chosen, as follows: 
\begin{equation*}
V_{1}\left( \mathbf{x}\right) =%
\begin{pmatrix}
v_{10}\left( \mathbf{x}\right) , & v_{11}\left( \mathbf{x}\right) , & \cdots
, & v_{1(N-1)}\left( \mathbf{x}\right)%
\end{pmatrix}%
^{T}\text{ for }v_{1n}=\hat{v}_{n}\left( \mathbf{x}\right) \tilde{\chi}%
\left( \mathbf{x}\right) ,
\end{equation*}%
where 
\begin{equation*}
\widetilde{\chi }\left( \mathbf{x}\right) =%
\begin{cases}
\exp \left( \frac{x^{2}}{x^{2}-b_{x}^{2}}\right) \exp \left( \frac{y^{2}}{%
y^{2}-b_{y}^{2}}\right) \exp \left( \frac{\left( z+R\right) ^{2}}{\left(
z+R\right) ^{2}-\left( R-b_{z}\right) ^{2}}\right) & \text{if }\mathbf{x}\in
\Omega _{2}, \\ 
0 & \text{otherwise}.%
\end{cases}%
\end{equation*}%
The second step allows us to obtain more accurate values of $c\left( \mathbf{%
x}\right) $ inside of the inclusion to be imaged. Our reconstruction is
concluded after we smooth the final solution $c_{\text{final}}\left( \mathbf{%
x}\right) $ by the Gaussian filtering via the \texttt{smooth3} function in 
\textsc{Matlab}. In particular, we find $c_{\text{comp}}\left( \mathbf{x}%
\right) $ as $c_{\text{comp}}\left( \mathbf{x}\right) =\text{smooth}%
(\left\vert c_{\text{final}}\left( \mathbf{x}\right) \right\vert (1+\widehat{%
p}))$, where $c_{\text{final}}\left( \mathbf{x}\right) $ is the function $%
c\left( \mathbf{x}\right) $ obtained in the last iterative step of the
minimization procedure of Step 2. This procedure is similar to the one of
smoothing $c_{\text{temp}}$ $\left( \mathbf{x}\right) $ on Step 1. Note that
due to (\ref{ccc}), $c_{\text{final}}\left( \mathbf{x}\right) >0$. Here,
finding the value of $\widehat{p}$ is based upon the maximal value of $%
\left\vert c_{\text{final}}\left( \mathbf{x}\right) \right\vert $. This
maximal value is computed with a good accuracy. Next, however, we need to
smooth the function $c_{\text{final}}\left( \mathbf{x}\right) $ using the
Gaussian filtering. The smoothed version of $c_{\text{final}}\left( \mathbf{x%
}\right) $ always has a lower maximal value. Therefore, we find such a
number $\widehat{p}\geq 0$ that $\max \left( c_{\text{temp}}\left( \mathbf{x}%
\right) \right) =\left( 1+\widehat{p}\right) \max \left( \text{smooth}\left(
c_{\text{temp}}\left( \mathbf{x}\right) \right) \right) $. Hence, the value
of $\widehat{p}$ varies in every single test. As an example, in Test 1 $\max
\left\vert c_{\text{final}}\left( \mathbf{x}\right) \right\vert =1.8873.$
But the smoothed version smooth$(\left\vert c_{\text{final}}\left( \mathbf{x}%
\right) \right\vert )$ without $\widehat{p}$ attains the maximal value
1.6446, while with $\widehat{p}=0.3765$ it moves back to the value $1.8873$.
The above forms the second step of the numerical solution of our CIP.

\textbf{Remark 7.1}. \emph{We point out that the computational procedure
described above does not use any advanced knowledge of a small neighborhood
of the correct solution. In other words, our reconstruction method converges 
\textbf{globally}; see\textbf{\ }Introduction for the definition of the
global convergence.}

We now mention that to run the minimization procedure, we need to compute
the gradient $J_{\lambda ,\gamma }^{^{\prime }}$ in (\ref{6.1}) of the
discrete functional $J_{\lambda ,\gamma }$ in \eqref{eq:discreteJ}. We have
discovered that having the expression for the gradient via an explicit
formula, significantly reduces the computational time. We have derived such
a formula using the technique of Kronecker deltas, which has been outlined
in \cite{Kuzhuget2010}. For brevity we do not provide this formula here.

\begin{table}[tbp]
\begin{centering}
		\begin{tabular}{|c|c|c|c|c|c|}
			\hline 
			Test & 1 & 2 & 3 & 4 & 5\tabularnewline
			\hline 
			Inclusion & Sphere & Ellipsoid & Ellipsoid & Rectangle & Sphere (noise)\tabularnewline
			\hline 
			$\max\left(c_{\text{true}}\right)$ & 2 & 5 & 10 & 2 & 2\tabularnewline
			\hline 
			$\max\left(c_{\text{comp}}\right)$ & 1.8873 & 5.1886 & 9.3461 & 2.1834 & 1.8782\tabularnewline
			\hline 
			Error (%
			\ref{7.2}) & 5.64\% & 3.77\% & 6.54\% & 9.17\% & 6.09\%\tabularnewline
			\hline 
		\end{tabular}
		\par\end{centering}
\caption{Description of numerical examples and the corresponding relative
errors between $\max\left(c_{\text{true}}\right)$ and $\max\left(c_{\text{%
comp}}\right)$. }
\label{table:1}
\end{table}
\begin{table}[tbp]
\begin{centering}
\begin{tabular}{|c|c|c|c|}
	\hline 
	Test & True center & Computed center & Lowest point\tabularnewline
	\hline 
	1 & $\left(0,0,-2.5\right)$ & $\left(0,0,-2.4\right)$ & $\left(0,0,-2.823\right)$\tabularnewline
	\hline 
	2 & $\left(0,0,-2\right)$ & $\left(-0.01,0.02,-1.98\right)$ & $\left(0,0,-2.81\right)$\tabularnewline
	\hline 
	3 & $\left(0,0,-2\right)$ & $\left(-0.02,-0.03,-2.13\right)$ & $\left(0,0,-2.58\right)$\tabularnewline
	\hline 
	4 & $\left(0,0,-2.5\right)$ & $\left(0,-0.01,-2.35\right)$ & $\left(-0.01,0,-2.79\right)$\tabularnewline
	\hline 
	5 & $\left(0,0,-2.5\right)$ & $\left(-0.01,0,-2.39\right)$ & $\left(0,0,-2.83\right)$\tabularnewline
	\hline 
\end{tabular}
\par\end{centering}
\caption{Comparison of the location information between the true and
computed objects. The true lowest point of all true objects is fixed at $%
\left(0,0,-2.8\right)$.}
\label{table:2}
\end{table}

\subsubsection{Reconstruction results}

\label{sub:7.2.3}

To this end, the dimensionless spatial variables are defined as $\mathbf{x}%
^{\prime }=\mathbf{x}/(10\text{ cm})$. This means that, for instance, in our
reference Test 1, where the mine-like target is ball-shaped with the radius $%
0.3,$ its actual diameter is $6$ cm. Therefore, with $R=3$ considered in %
\cref{sec:experiments} we suppose to look for a possible explosive in a
cubic area $0.216\text{ m}^{3}$, where $m$ stands for meters. We set in (\ref%
{2001}) $d=7.5$, which means that the distance between the line of sources,
whose length is $20$ cm, and the ground surface should be $45$ cm, which is
realistic for some ground penetrating radars.


Now, since $k=2\pi /\upsilon ,$ where $\upsilon $ is the wave length, then,
after the change of variables $\mathbf{x}^{\prime }=\mathbf{x}/(10$ cm$)$ in
the Helmholtz equation \eqref{eq:forward1}, the dimensionless number $k=2\pi
\cdot \left( 10\text{ cm}\right) /\upsilon .$ Following the previous
publication of this group about the work with experimental data \cite%
{Nguyen2018}, we use here the frequency of 3.15 GHz, which means that the
wavelength $\upsilon =9.5$ cm. Therefore, the dimensionless value of $k$ we
work with is 
\begin{equation*}
k=2\pi \cdot \left( 10\text{ cm}\right) /\left( 9.5\text{ cm}\right) =6.6.
\end{equation*}

In our Tests 1--5 listed below, we want to accurately reconstruct all three
components of targets: locations, target/background contrasts and shapes.

In all our tests, the mine-like objects are such that the lowest point of
their front surfaces is fixed at $(x,y,z)=(0,0,-2.8).$ However, when running
the minimization procedure, we do not assume a knowledge of neither the
locations nor the shapes of inclusions. We only assume that any inclusion of
our interest is located close to the ground surface $\Gamma $, which
concerns the choice of the smooth function $\chi $ in \cref{sec:experiments}.

It was demonstrated in \cite{Nguyen2018} that the data propagation procedure
mentioned in \cref{sec:3.3} estimates quite accurately the $z-$coordinate of
the front surface of the target. Hence, we indeed can assume that the plane $%
\left\{ z=-R\right\} $ is close to the target to be imaged. On the other
hand, if it would be far, then the image would be worse. Indeed, in this
case the backscattering data at $\Gamma $ would be less sensitive to the
presence of the inclusion and, at the same time, the transmitted data at $%
\overline{\Omega }\cap \left\{ z=R\right\} $ are artificially set not to be
sensitive to this presence.


In the reconstruction results, we are concerned with the relative error
between $\max \left( c_{\text{true}}\left( \mathbf{x}\right) \right) $ and $%
\max \left( c_{\text{comp}}\left( \mathbf{x}\right) \right) $, where $c_{%
\text{comp}}\left( \mathbf{x}\right) $ is the computed function $c\left( 
\mathbf{x}\right) $. More precisely, we define this error as 
\begin{equation}
E_{\max }=\frac{\left\vert \max \left( c_{\text{true}}\left( \mathbf{x}%
\right) \right) -\max \left( c_{\text{comp}}\left( \mathbf{x}\right) \right)
\right\vert }{\max \left( c_{\text{true}}\left( \mathbf{x}\right) \right) }%
\times 100\%.  \label{7.2}
\end{equation}%
Values of $\max \left( c_{\text{true}}\left( \mathbf{x}\right) \right) ,\max
\left( c_{\text{comp}}\left( \mathbf{x}\right) \right) $ and $E_{\max }$ for
all five tests are tabulated in \cref{table:1}. In the following, we depict
both three-dimensional true and computed inclusions by using the \texttt{%
isosurface} function in \textsc{Matlab} with the associated \texttt{isovalue}
being 5\% of the maximal value. As to the locations of computed inclusions,
we briefly report them in \cref{table:2}.

\begin{figure}[]
\centering
\begin{subfigure}[b]{0.32\linewidth}
		\includegraphics[width=\linewidth]{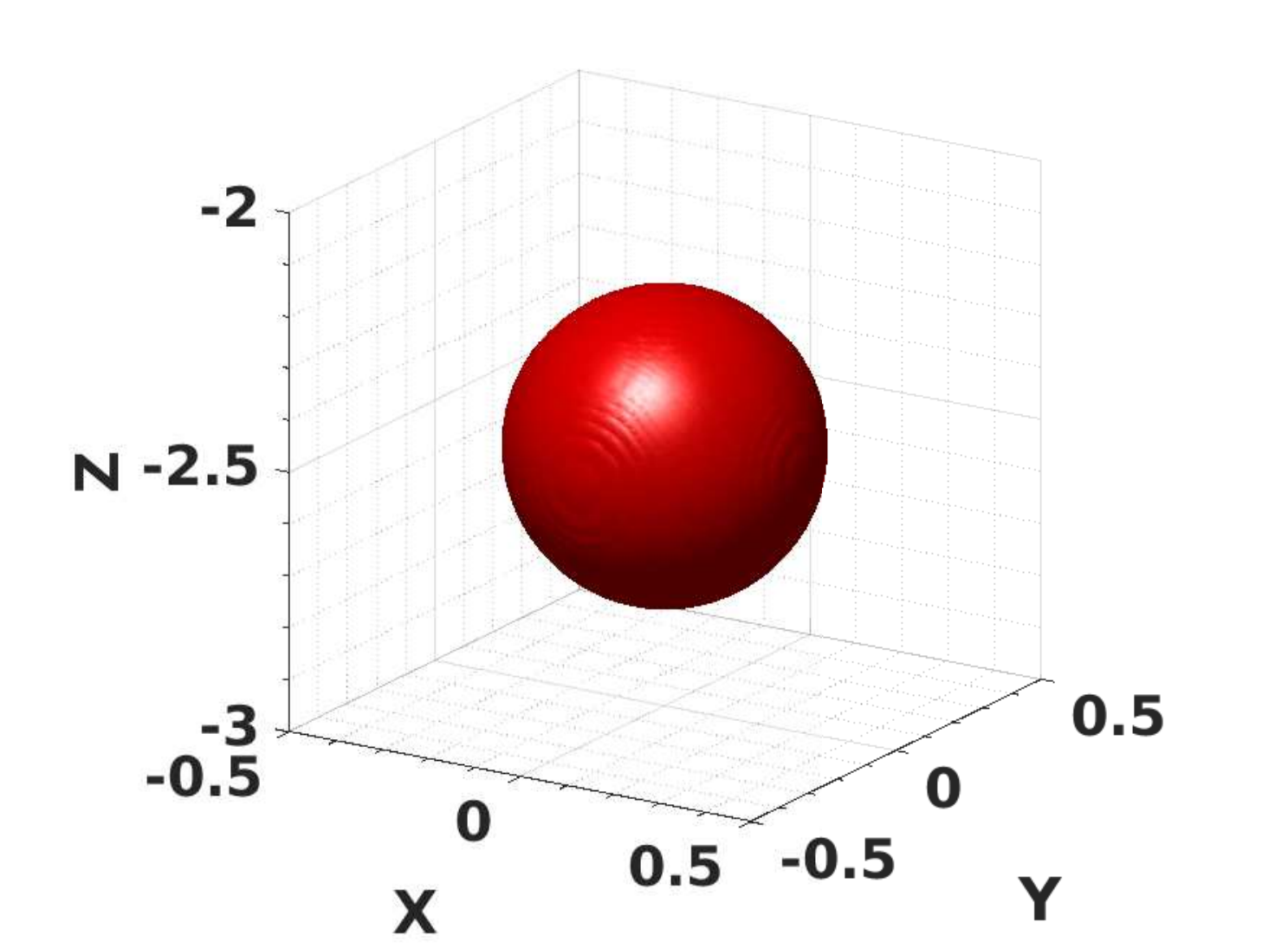}
		\caption{$c_{\text{true}}=2$\label{fig:22}}
	\end{subfigure}
\begin{subfigure}[b]{0.32\linewidth}
		\includegraphics[width=\linewidth]{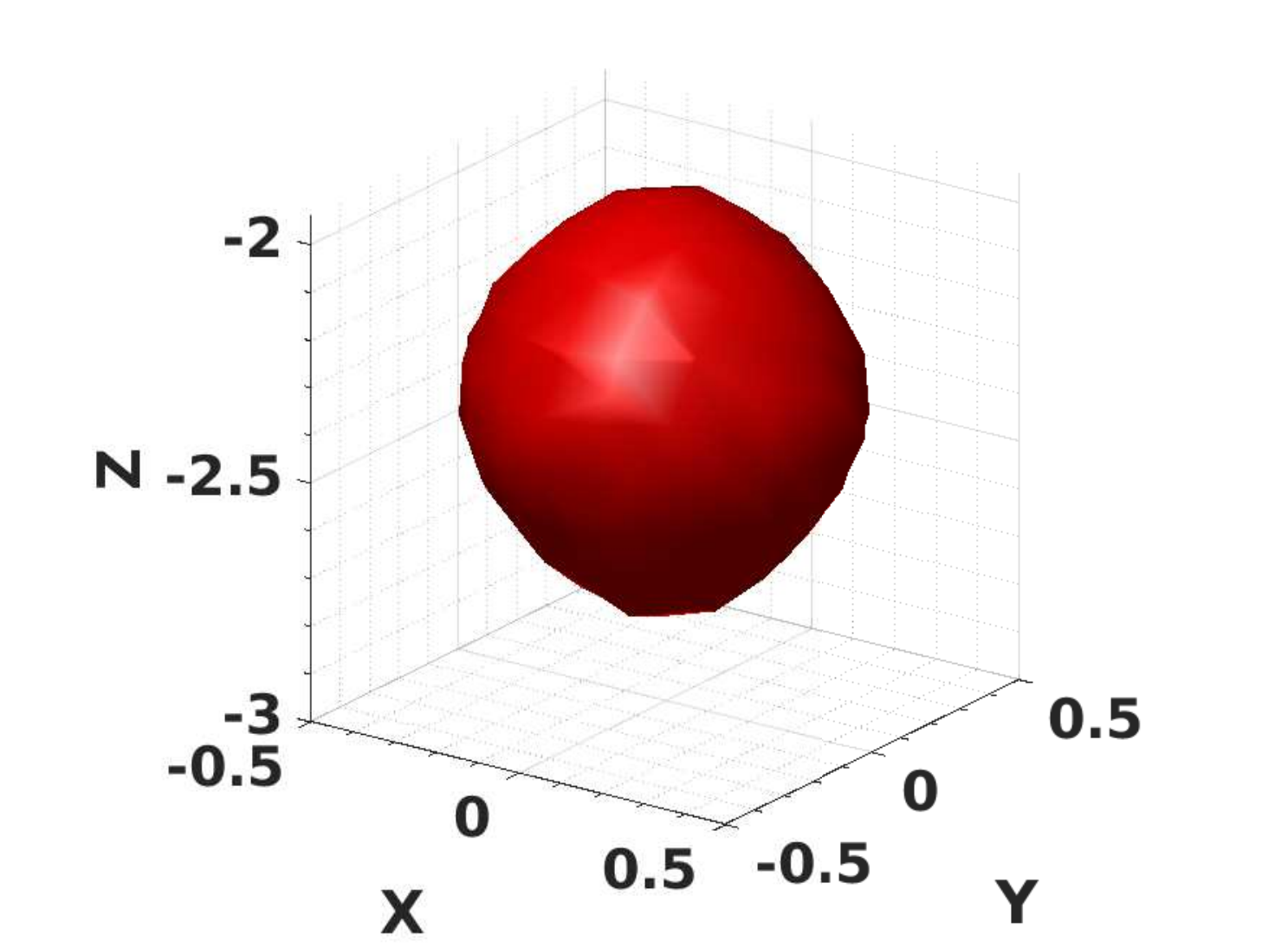}
		\caption{$c_{\text{comp}}=1.8873$\label{fig:1a}}
	\end{subfigure}
\begin{subfigure}[b]{0.32\linewidth}
		\includegraphics[width=\linewidth]{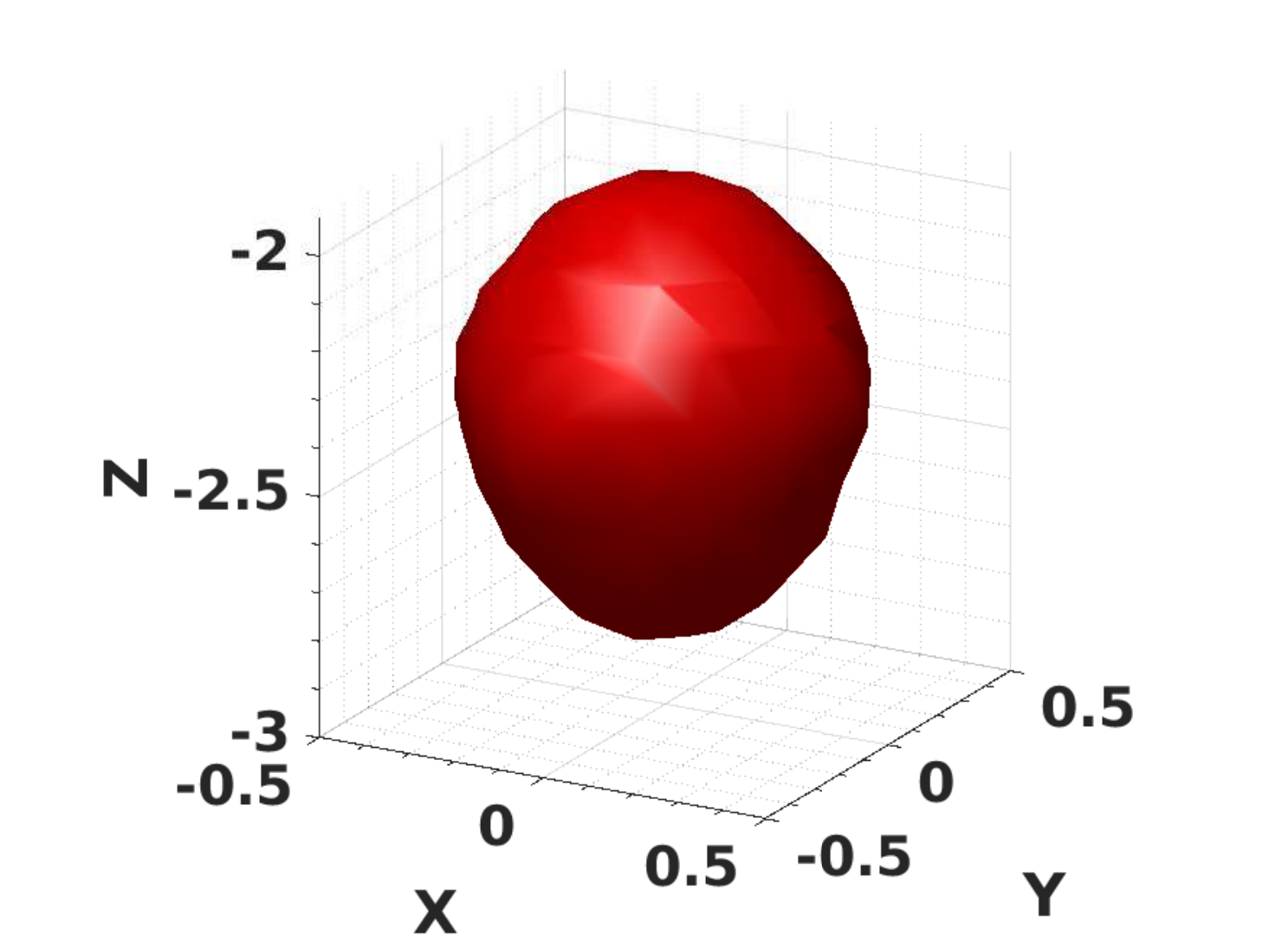}
		\caption{$c_{\text{comp}}=1.8782$\label{fig:5bb}}
	\end{subfigure}
\caption{Reconstructions of a ball-shaped object with the correct value of
the coefficient $c_{\text{true}}=2$. Without noise (Test 1), the maximal
value $\max \left( c_{\text{comp}}\right) $ of the computed coefficient $c$
is 1.8873. With 5$\%$ noise (Test 5), the maximal value $\max \left( c_{%
\text{comp}}\right) $ of the computed coefficient $c$ is 1.8782. (a) True
image. (b) Reconstructed image of Test 1. (c) Reconstructed image of Test 5
with noisy data.}
\label{fig:3}
\end{figure}


\subsubsection*{Test 1. Ball-shaped inclusion}

In this test, we examine our numerical method for the case of a ball-shaped
object with the dielectric constant in it $c_{\text{true}}=2$. Images of the
true object and its reconstruction are presented in \cref{fig:3}. We observe
in \cref{fig:22} that the shape of the reconstructed object is imaged
accurately. We also obtain that the lowest point of the computed one is at $%
(0,0,-2.823)$, while it should be $(0,0,-2.8) $ for the true solution as we
have set up above. So, these are close. The size of the computed inclusion
is slightly larger than the correct one. Together with the approximate
dielectric constant, we conclude that location, shape and $\max \left( c_{%
\text{comp}}\right) $ are reconstructed accurately.

\subsubsection*{Tests 2--3. Ellipsoids}

\begin{figure}[]
\centering
\begin{subfigure}[b]{0.32\linewidth}
		\includegraphics[width=\linewidth]{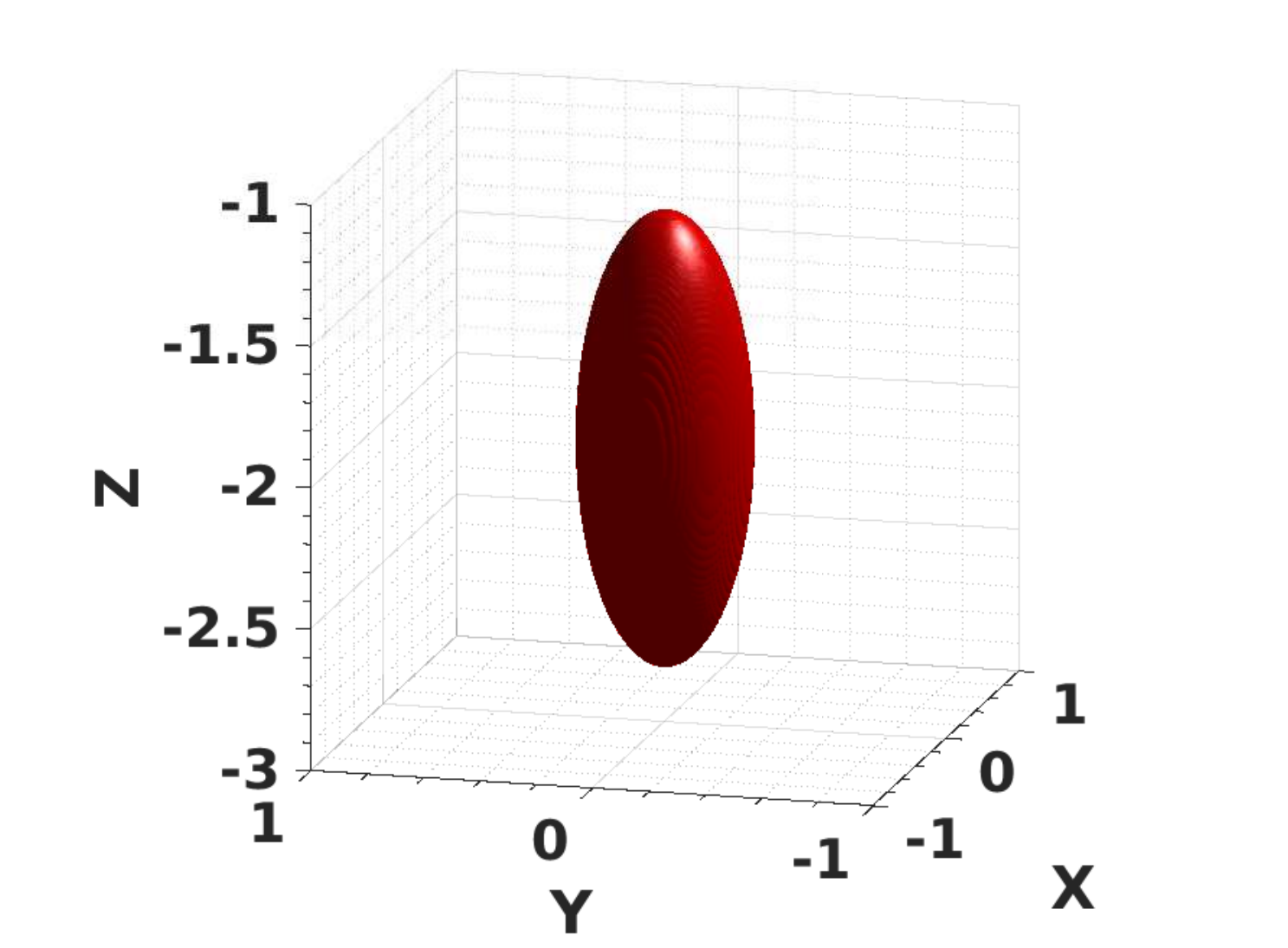}
		\caption{$c_{\text{true}}=5,10$}
	\end{subfigure}
\begin{subfigure}[b]{0.32\linewidth}
		\includegraphics[width=\linewidth]{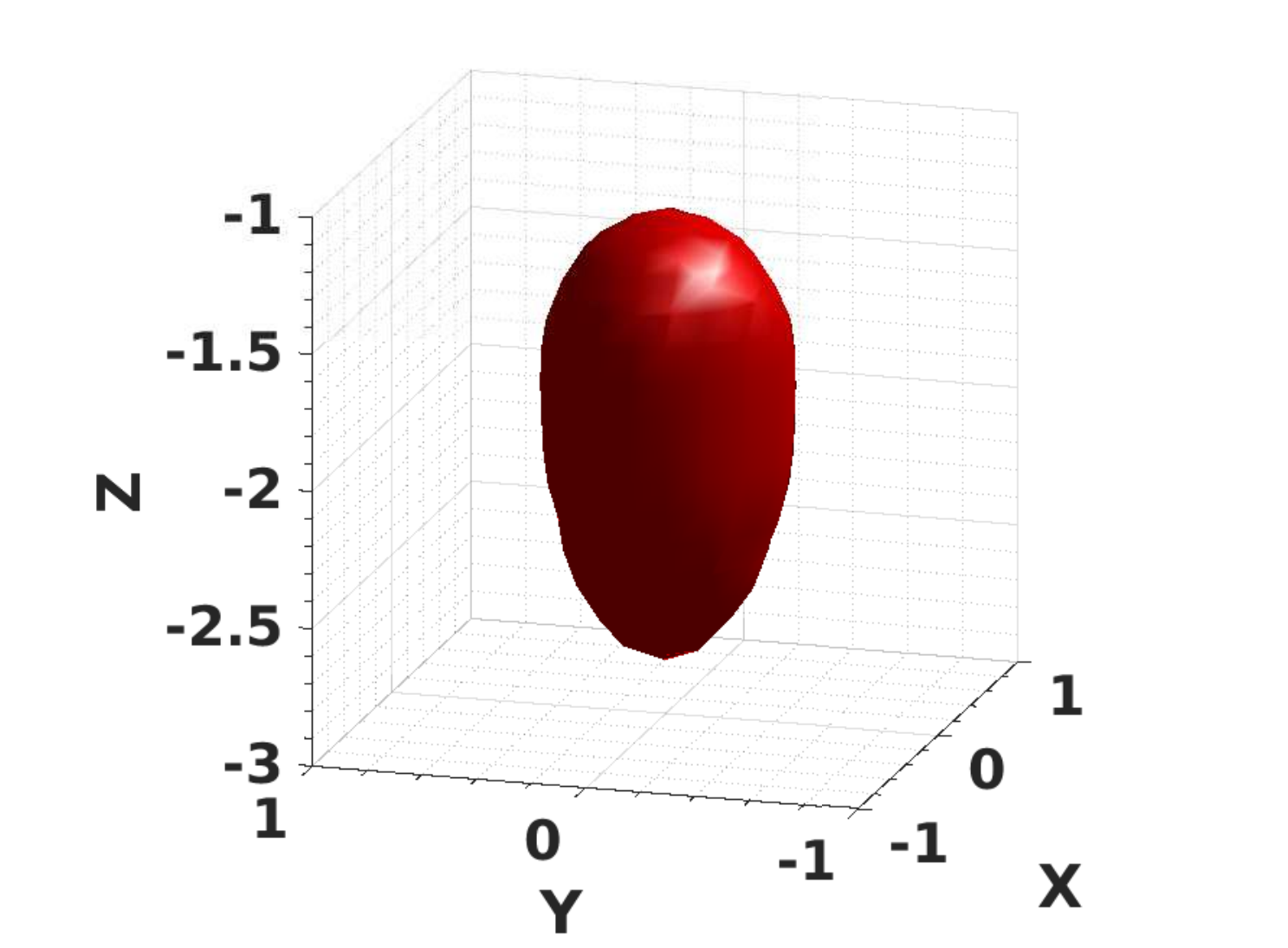}
		\caption{$c_{\text{comp}}=5.1886$}
	\end{subfigure}
\begin{subfigure}[b]{0.32\linewidth}
		\includegraphics[width=\linewidth]{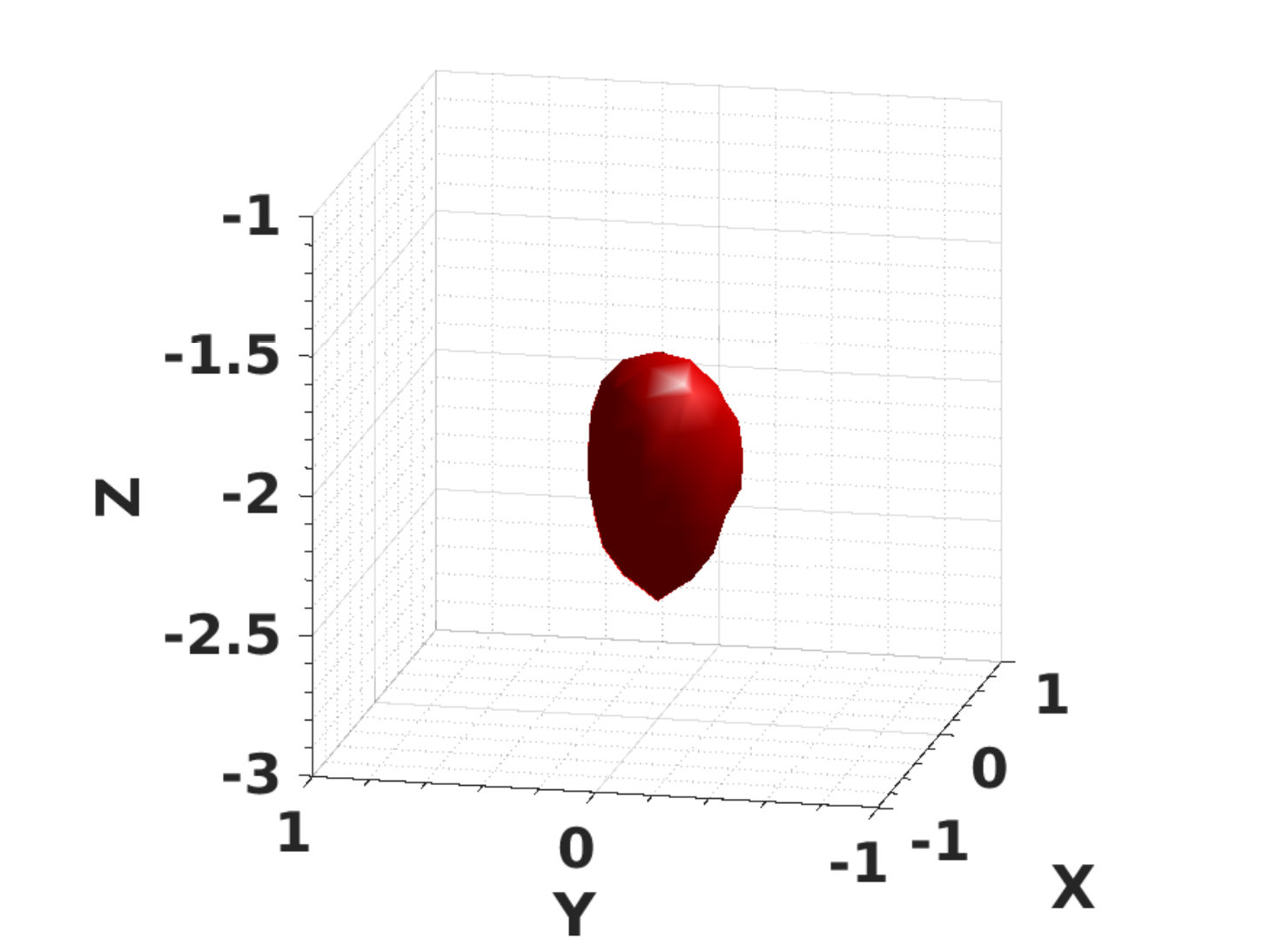}
		\caption{$c_{\text{comp}}=9.3461$\label{fig:5a}}
	\end{subfigure}
\caption{Reconstructions of an ellipsoid with the correct values of the
coefficient $c_{\text{true}}=5,10$ (Tests 2 and 3). The maximal values $\max
\left( c_{ \text{comp}}\right) $ of the computed coefficient $c$ are 5.1886
and 9.3461, respectively. (a) True image. (b) Reconstructed image of Test 2.
(c) Reconstructed image of Test 3. See the text for an explanation of the
decrease of the size of image (c), compared with (b).}
\label{fig:4}
\end{figure}

The next two tests, Tests 2 and 3, are about ellipsoidal targets with high
target/background contrast levels of 5 and 10. Here, we consider the
ellipsoids with principal semi-major axis and two semi-minor axes,
respectively, being 0.8 and 0.3. The objects are centered at $(0,0,-2)$. The
computational results show that the sizes of reconstructed ellipsoidal
targets with $c_{\text{true}}=5,10$ decrease when $c_{\text{true}}$
increases. Here is an explanation on this. We have computationally observed
that the area of the region on $\Gamma $, where the data concentrate around
the maximal absolute value, decreases when $c_{\text{true}}$ increases from
5 to 10. We still have observed the ellipsoidal shape of the reconstructed
target; see \cref{fig:5a}. Besides, the center of the computed one with $c_{%
\text{true}}=10$ is close to the true position and the number $\max \left(
c_{\text{comp}}\left( \mathbf{x}\right) \right) $ is still good; see, e.g. %
\cref{table:1,table:2}. Therefore, the results for $c_{\text{true}}=10$
still accommodate the above listed main purposes of our reconstructions.
Last but not least, the reconstruction of the number $c_{\text{true}}=5$ is
very accurate.



\textbf{Remark 7.2}. \emph{We point out that inclusions with high contrasts
of, e.g., Tests 2--3 are unlikely to be imaged by conventional techniques.}

\subsubsection*{Test 4. Rectangular prism}

\begin{figure}[htbp]
\centering
\begin{subfigure}[b]{0.32\linewidth}
		\includegraphics[width=\linewidth]{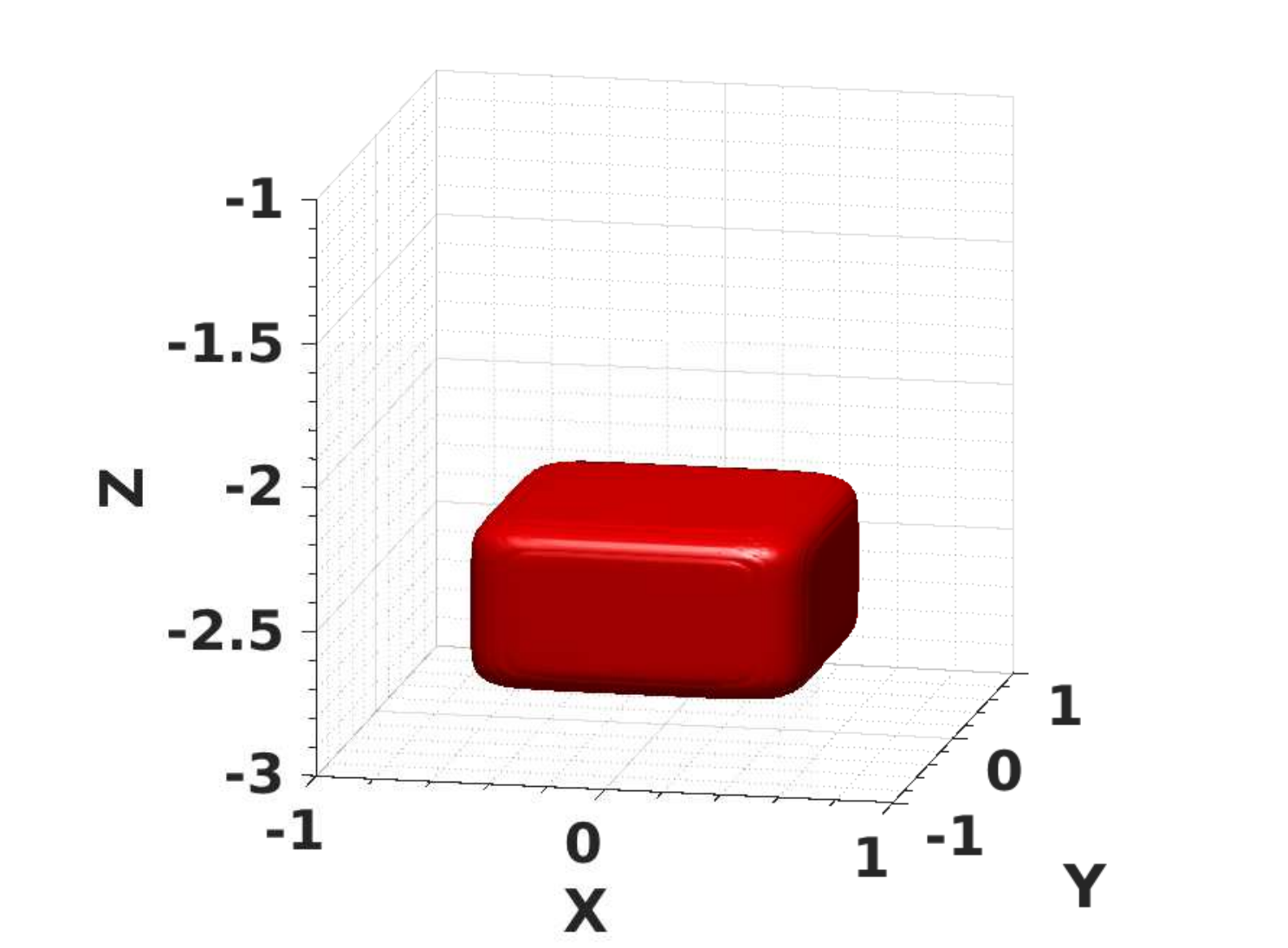}
		\caption{$c_{\text{true}}=2$}
	\end{subfigure}
\begin{subfigure}[b]{0.32\linewidth}
		\includegraphics[width=\linewidth]{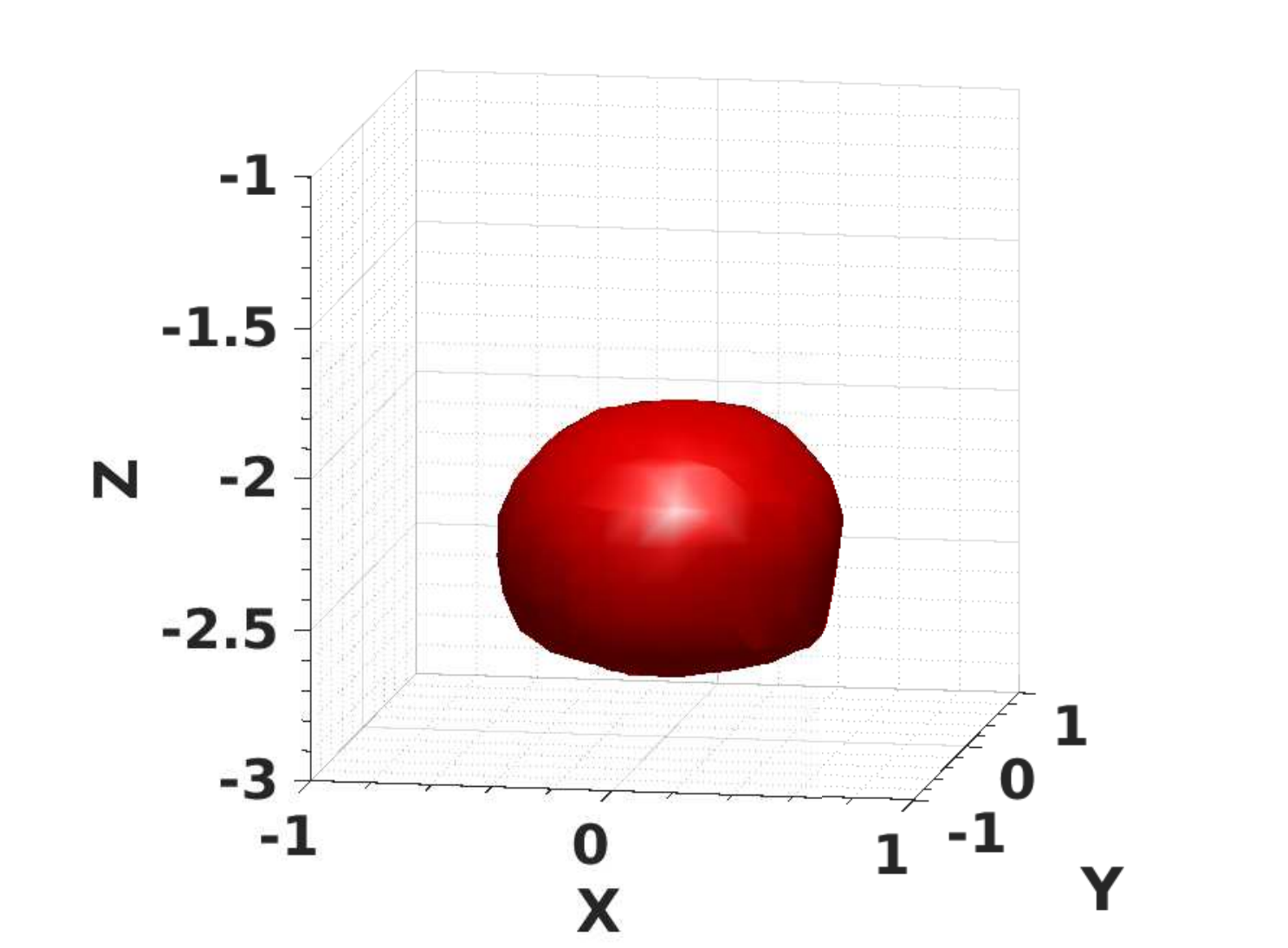}
		\caption{$c_{\text{comp}}=2.1834$}
	\end{subfigure}
\caption{Reconstruction of a rectangular prism with the correct value of the
coefficient $c_{\text{true}}=2$ (Test 4). The maximal value $\max \left( c_{ 
\text{comp}}\right) $ of the computed coefficient $c$ is 2.1834. (a) True
image. (b) Reconstructed image.}
\label{fig:6}
\end{figure}

In this test, we consider a rectangular prism with dielectric constant $c_{%
\text{true}}=2$. This object is of the size $1.2\times 1.2\times 0.6$ (W$%
\times $L$\times $H), which is realistic for the above-mentioned targets.
The reconstruction result for this case is displayed on \cref{fig:6}. The
lowest point of the reconstructed object is located at the point $%
(0,0,-2.79) $ while the correct point is $(0,0,-2.8).$ We see that this
result still fulfills the main purpose of our reconstructions listed above
in this section; see \cref{table:1,table:2} and \cref{fig:6}. It should be
noticed that this prism's boundary regularity is mathematically weaker than
the one of the spherical object. This explains why we obtain a smooth shape
rather than the one with a sharp boundary.

\subsubsection*{Test 5. Ball-shaped inclusion with a noise in the data}

In the last Test 5, we concentrate on the reconstruction with noisy data $F$
and $G$ defined, respectively, in (\ref{303}) and \cref{sec:3.3}. In doing
so, we simply add a random multiplicative noise to the simulated Dirichlet
and Neumann data obtained when solving the Lippmann--Schwinger equation: 
\begin{equation*}
F_{\text{noise}}\left( \mathbf{x},\alpha \right) =F\left( \mathbf{x},\alpha
\right) \left( 1+\delta \text{rand}\right) ,\quad G_{\text{noise}}\left( 
\mathbf{x},\alpha \right) =G\left( \mathbf{x},\alpha \right) \left( 1+\delta 
\text{rand}\right) .
\end{equation*}%
Here, $\delta \in (0,1)$ represents the noise level and \textquotedblleft $%
\text{rand"}$ is a random number uniformly distributed in the interval $%
(-1,1)$. Thus, this addition is pointwise with respect to both spatial
points and the point source. Given $\delta =0.05,$ which corresponds to the $%
5\%$ noise, we choose to use this noisy data in our Test 1. One can see from %
\cref{fig:5bb} that our reconstruction still has a good performance.

%
%
%




\end{document}